\documentclass{article}

\usepackage{amssymb}
\usepackage{amsmath}
\usepackage[all]{xypic}
\usepackage[frenchb]{babel}
\usepackage{a4wide}
\usepackage{enumerate, fancyhdr}

\begin{document}

\title{Torsion dans un produit de courbes elliptiques}

\author{Marc Hindry\footnote{hindry@math.jussieu.fr}, Nicolas Ratazzi \footnote{nicolas.ratazzi@math.u-psud.fr}
}

\renewcommand{\epsilon}{\varepsilon}
\renewcommand{\tilde}{\widetilde}

\newcommand{\Gal}{\textnormal{Gal}}
\newcommand{\tors}{\textnormal{tors}}
\newcommand{\Ker}{\textnormal{Ker}}
\newcommand{\red}{\textnormal{r\'{e}d}}
\newcommand{\Z}{\mathbb{Z}}
\newcommand{\Q}{\mathbb{Q}}
\newcommand{\C}{\mathbb{C}}
\newcommand{\R}{\mathbb{R}}
\newcommand{\N}{\mathbb{N}}
\renewcommand{\P}{\mathbb{P}}
\newcommand{\F}{\mathbb{F}}
\newcommand{\G}{\mathbb{G}}
\newcommand{\SD}{\mathbb{S}}
\renewcommand{\hat}{\widehat}

\newcommand{\T}{\mathcal{T}}

\newcommand{\K}{\overline{K}}

\newcommand{\mP}{\mathcal{P}}
\newcommand{\mA}{\mathcal{A}}
\renewcommand{\O}{\mathcal{O}}
\renewcommand{\L}{\mathcal{L}}
\newcommand{\V}{\mathcal{V}}
\newcommand{\kk}{\mathbb{K}}

\newcommand{\im}{\textnormal{Im}}

\newcommand{\m}{\mathfrak{m}}
\newcommand{\I}{\mathfrak{I}}

\newcommand{\VV}{\mathsf{V}}

\newcommand{\mmid}{\|\|}

\newcommand{\Proj}{\textnormal{Proj\,}}
\newcommand{\Spec}{\textnormal{Spec\,}}
\newcommand{\Deg}{\textnormal{Deg\,}}
\newcommand{\e}{\varepsilon}
\newcommand{\ab}{\textnormal{ab}}
\newcommand{\End}{\textnormal{End}}
\renewcommand{\H}{\textnormal{H}}
\newcommand{\MT}{\textnormal{MT}}
\newcommand{\codim}{\textnormal{codim\,}}
\newcommand{\Hom}{\textnormal{Hom}}
\newcommand{\Aut}{\textnormal{Aut}}
\newcommand{\Hdg}{\textnormal{Hdg}}
\newcommand{\GL}{\textnormal{GL}}
\newcommand{\SL}{\textnormal{SL}}
\newcommand{\Sp}{\textnormal{Sp}}
\newcommand{\GSp}{\textnormal{GSp}}
\newcommand{\SLG}{\textnormal{SLG}}
\newcommand{\Res}{\textnormal{Res}}
\renewcommand{\ss}{\sigma}
\newcommand{\pr}{\textnormal{pr}}
\newcommand{\hdg}{\textnormal{Hdg}}
\newcommand{\Kc}{K^{\text{cycl}}}
\newcommand{\ssl}{\mathfrak{sl}}

\renewcommand{\text}{\textnormal}

\newtheorem{theo}{Th{\'e}or{\`e}me} [section]
\newtheorem{lemme}[theo]{Lemme}
\newtheorem{conj}{Conjecture}[section]
\newtheorem{prop}[theo]{Proposition}
\newtheorem{cor}[theo]{Corollaire}
\newtheorem{quest}[theo]{Question}

\newcommand{\defi}{\addtocounter{theo}{1}{\noindent \textbf{D{\'e}finition \thetheo\ }}}
\newcommand{\rem}{\addtocounter{theo}{1}{\noindent \textbf{Remarque \thetheo\ }}}

\newcommand{\demo}{\noindent \textit{D{\'e}monstration} : }
\newcommand{\findemo}{\hfill \Box}

\maketitle

\hrulefill

\bigskip

\noindent \textbf{R\'esum\'e : } 
Soit $A$ une vari\'et\'e ab\'elienne d\'efinie sur un corps de nombres $K$, le nombre de points de torsion
d\'efinis sur une extension finie $L$   est born\'e polynomialement en terme du degr\'e $[L:K]$. Nous formulons une question sug\'erant l'exposant optimal  dans cette borne en terme de la dimension du groupe de Mumford-Tate des sous-vari\'et\'es ab\'eliennes de $A$; nous \'etudions le comportement par produit  et r\'epondons par l'affirmative \`a la question dans le cas d'un produit de courbes elliptiques.
  
\bigskip

\hrulefill

\bigskip

\noindent \textbf{Abstract : }
 Let $A$ be an abelian variety defined over a number field $K$, the number of torsion points rational over a finite extension $L$ is bounded polynomially in terms of the degree $[L:K]$. We formulate a question suggesting the optimal exponent for this bound in terms of the dimension of the Mumford-Tate groups of the abelian subvarieties of $A$; we study the behaviour under product and then give a positive answer to our question when $A$ is the product of elliptic curves.

\hrulefill

\section{Introduction et r\'esultats}

Soit $A/\overline{\Q}$ une vari\'et\'e ab\'elienne, d\'efinie sur un corps de nombres $K$, de dimension $g\geq 1$. Le classique th\'eor\`eme de Mordell-Weil assure que le groupe $A(K)$ des points $K$-rationnels de $A$ est de type fini. Un probl\`eme naturel qui se pose alors est de comprendre le sous-groupe de torsion $A(K)_{\tors}$. Un premier probl\`eme consiste en fait \`a essayer de bien appr\'ehender le cardinal de $A(K)_{\tors}$ lorsque $A$ et/ou $K$ varient. Comme dans l'article \cite{ratazziens} auquel ce papier fait suite, nous nous int\'eressons ici au cas o\`u l'on fixe une vari\'et\'e ab\'elienne $A$ d\'efinie sur un corps de nombres $K_0$ et o\`u l'on fait varier $K$ parmi les extensions finies de $K_0$ ; l'objectif \'etant cette fois-ci d'obtenir une borne  avec une d\'ependance explicite et, si possible, optimale  en le degr\'e $[K:K_0]$ ou, ce qui revient au m\^eme, en le degr\'e $[K:\Q]$. Concernant cette question, Masser \cite{lettre} et \cite{mas} a montr\'e dans le cas g\'en\'eral que la d\'ependance est polynomiale en $[K:\Q]$; plus pr\'ecis\'ement il montre l'\'enonc\'e suivant.

\begin{theo}\label{theomasser} \textnormal{\textbf{(Masser  \cite{lettre})}} Soit $A$ une vari\'et\'e ab\'elienne de dimension $g$, d\'efinie sur un corps de nombres $K_0$, il existe une constante $c_A$ telle que, pour toute extension finie $K$ de $K_0$ on a~:
$$\left|A(K)_{\tors}\right|\leq c_A[K:\Q]^g(\log [K:\Q])^g.$$
\end{theo}
Masser indique d'ailleurs que l'exposant $g$ n'est probablement pas le meilleur possible (sauf pour le cas d'une puissance d'une courbe elliptique \`a multiplication complexe).
 La question naturelle qui se pose est alors de savoir quel est le plus petit exposant $\gamma(A)$ possible dans cette borne polynomiale. Dans \cite{ratazziens}, le second auteur a donn\'e une r\'eponse \`a cette question dans le cas des vari\'et\'es ab\'eliennes CM sans facteur carr\'e (\textit{i.e.} de la forme $\prod A_i$, les $A_i/\overline{\Q}$ \'etant de type CM, deux \`a deux non-isog\`enes). 

\medskip

\defi Soit $A/\overline{\Q}$ une vari\'et\'e ab\'elienne d\'efinie sur un corps de nombres $K_0$. On pose 
\[\gamma(A)=\inf\left\lbrace x>0\, | \,  \forall F/K_0 \text{ finie, }\ \left|A(F)_{\tors}\right|\ll [F:K_0]^x\right\rbrace.\]

\noindent  La notation $\ll$ signifie qu'il existe  une constante  $C$, ne d\'ependant que de $A/K_0$, telle que l'on a  $\left|A(F)_{\tors}\right|\leq C [F:K_0]^x$. On voit facilement que l'invariant d\'efini ci-dessus est ind\'ependant du corps de d\'efinition $K_0$ choisi et ne d\'epend en fait que de la classe d'isog\'enie de la vari\'et\'e ab\'elienne $A$.

\noindent Introduisons par ailleurs un autre invariant, d\'efini en terme de la dimension du groupe de Mumford-Tate d'une vari\'et\'e ab\'elienne (la d\'efinition de ce groupe est rappel\'ee au paragraphe \ref{mtt}).

\medskip

\noindent  Soit $A/\overline{\Q}$ une vari\'et\'e ab\'elienne isog\`ene au produit $\prod_{i=1}^n A_i^{n_i}$ o\`u les $A_i$ sont des vari\'et\'es ab\'eliennes simples deux \`a deux non-isog\`enes et o\`u les $n_i$ sont des entiers strictement positifs. 

\medskip

\defi On d\'efinit l'invariant $\alpha(A)$ par
\[\alpha(A)=\max_{\emptyset\not= I\subset\{1,\ldots,n\}}\frac{2\sum_{i\in I}n_i\dim A_i}{\dim\MT\left(\prod_{i\in I}A_i\right)}.\]

\medskip

\begin{quest}\label{conjR} Soit $A/K$ est une vari\'et\'e ab\'elienne sur un corps de nombres. A-t-on $\alpha(A)=\gamma(A)$~?
\end{quest}

\noindent On remarquera (\textit{cf.} la proposition \ref{mtmin} plus bas) que, sauf dans le cas d'une puissance d'une courbe elliptique \`a multiplication complexe (\textit{i.e.} $A=A_1^n$ avec $\dim A_1=1$ et $\dim\MT(A)=\dim\MT(A_1)=2$), l'invariant $\alpha(A)$ est strictement plus petit que $\dim A$.

\medskip

\noindent L'objet de la suite de cet article est de donner une r\'eponse affirmative \`a la question \ref{conjR} ci-dessus  si $A$ est un produit de courbes elliptiques. Commen\c{c}ons par faire une simple remarque (dont la preuve est donn\'ee au paragraphe \ref{pp1}).

\begin{prop}\label{p1} Si $A/\overline{\Q}$ est une vari\'et\'e ab\'elienne quelconque, l'in\'egalit\'e $\gamma(A)\geq\alpha(A)$ est vraie.
\end{prop}

\noindent Ceci \'etant nous pouvons \'enoncer notre r\'esultat principal~:

\begin{theo}\label{th1} Si $A=\prod_{i=1}^mE_i$ est un produit de courbes elliptiques sur $\overline{\Q}$, isog\`enes ou non, alors
\[\gamma(A)=\alpha(A).\]
De plus, si les $E_i$ sont deux \`a deux non-isog\`enes, ce nombre vaut $2m/(1+3m)$ si toutes les $E_i$ sont sans multiplication complexe et
$2r/(1+r)$ si $E_1,\dots,E_r$ sont \`a multiplication complexe et $E_{r+1},\dots,E_{m}$ sont sans multiplication complexe.
\end{theo}

\medskip

\noindent Notons que dans le cas d'un produit de courbes elliptiques sans facteur carr\'e (\textit{i.e.} $A$ isog\`ene \`a un produit $\prod_{i=1}^mE_i$, les $E_i$ \'etant deux \`a deux non-isog\`enes) ce th\'eor\`eme nous donne notamment la majoration 
\[\gamma(A)\leq\begin{cases}
\frac{2}{3}	&	\text{ si $A$ est sans CM}\\
2			&	\text{sinon},
\end{cases}\]
\noindent alors que le r\'esultat de Masser donnait simplement $\gamma(A)\leq \dim A$.

\medskip 

\noindent Avant de d\'emontrer le th\'eor\`eme \ref{th1}, nous faisons un certain nombre de remarques et pr\'elimi\-nai\-res que nous formulons dans le contexte g\'en\'eral des vari\'et\'es ab\'eliennes lorsque cela est possible. La partie 2 est consacr\'ee \`a divers rappels et r\'esultats concernant le groupe de Mumford-Tate et la conjecture de Mumford-Tate. Le paragraphe 4 explique quelques r\'eductions du probl\`eme : on y montre notamment que, si $A=\prod_{i=1}^rA_i^{n_i}$, il suffit d'obtenir, pour $H$ sous-groupe fini produit de $A[\ell^{\infty}]$, une borne du cardinal de $H$ en fonction du degr\'e de l'extension $[K(H):K]$. Les paragraphes 5 et 6 sont consacr\'es \`a des arguments galoisiens permettant d'entamer le d\'evissage du calcul de $\gamma(A)$, avec $A=\prod_{i=1}^rA_i^{n_i}$. Nous indiquons une m\'ethode permettant, sous certaines hypoth\`eses sur $A$, de ramener le probl\`eme \`a une combinatoire ne faisant intervenir que la valeur des $\gamma(A_i)$ pour les facteurs simples $A_i$ de $A$. Jusqu'au paragraphe 6 compris, les arguments sont valables dans le cadre g\'en\'eral d'un produit de vari\'et\'es  ab\'eliennes. Au paragraphe 7, nous expliquons comment achever le calcul de $\gamma(A)$ pour un produit de courbes elliptiques, utilisant la technique d\'evelopp\'ee dans les paragraphes pr\'ec\'edents, ainsi que les r\'esultats d\'ej\`a connus pour les courbes elliptiques.

\section{Rappels sur les groupes et la conjecture de Mumford-Tate\label{mtt}}

\noindent Fixons d\'esormais un corps de nombres $K$ et un plongement $K\subset \C$ et notons $\overline{K}$ une cl\^oture alg\'ebrique de $K$ dans $\C$. Soit $A/K$ une vari\'et\'e ab\'elienne. On note $V=H^1(A(\C),\Q)$ le premier groupe de cohomologie singuli\`ere de la vari\'et\'e analytique complexe $A(\C)$. C'est un $\Q$-espace vectoriel de dimension $2g$. Il est naturellement muni d'une structure de Hodge de type $\{(1,0),(0,1)\}$, c'est-\`a-dire d'une d\'ecomposition sur $\C$ de $V_{\C}:=V\otimes_{\Q}\C$ donn\'ee par $V_{\C}=V^{1,0}\oplus V^{0,1}$ telle que $V^{0,1}=\overline{V^{1,0}}$ o\`u $\overline{\hspace{.1cm} \cdot \hspace{.1cm}}$ d\'esigne la conjugaison complexe. On note $\mu : \G_{m,\C}\rightarrow \GL_{V_{\C}}$ le cocaract\`ere tel que pour tout $z\in\C^{\times}$, $\mu(z)$ agit par multiplication par $z$ sur $V^{1,0}$ et agit trivialement sur $V^{0,1}$. On d\'efinit le groupe de Mumford-Tate en suivant \cite{pink1}.

\medskip

\defi Le \textit{groupe de Mumford-Tate} $\MT(A)/\Q$ de $A$ est le plus petit $\Q$-sous-groupe alg\'ebrique $G$ de $\GL_V$ (vu comme $\Q$-sch\'ema en groupes) tel que, apr\`es extension des scalaires \`a $\C$, le cocaract\`ere $\mu$ se factorise \`a travers $G_{\C}:=G\times_{\Q}\C$.

\medskip

\noindent  Concernant les groupes de Mumford-Tate $\MT(A)$ l'\'enonc\'e suivant est bien connu mais nous n'avons pas trouv\'e de r\'ef\'erence.

\begin{lemme}\label{puissance0} Soit $\prod_{i=1}^rA_i^{n_i}/K$ une vari\'et\'e ab\'elienne sur un corps de nombres $K$ plong\'e dans $\C$. On a 
\[\MT\left(\prod_{i=1}^rA_i^{n_i}\right)\simeq\MT\left(\prod_{i=1}^rA_i\right).\]
\end{lemme}
\demo Remarquons par la formule de K\"unneth, que 
\[H^1\left(\prod_{i=1}^rA_i^{n_i}(\C),\Q\right)\simeq \bigoplus_{i=1}^r H^1\left(A_i(\C),\Q\right)^{\oplus n_i}.\]
\noindent Posons $V_i=H^1(A_i(\C),\Q)$. Vu la d\'efinition du groupe de Mumford-Tate, on voit que l'on obtient le r\'esultat en plongeant diagonalement $\GL_{V_i}$ dans $\GL_{V_i^{\oplus n_i}}$ et $\GL_{V_i^{\oplus n_i}}$ dans $\GL_{V}$ avec $V=\bigoplus_{i=1}^{r}V_i^{\oplus n_i}$.\hfill$\Box$

\medskip

\noindent Nous rappelons \'egalement un r\'esultat concernant le lien entre groupe de Mumford-Tate du produit et produit des groupes de Mumford-Tate, dans le cas des courbes elliptiques. Pr\'ecis\'ement, dans ce cas, le bon objet n'est pas le groupe de Mumford-Tate, mais le groupe de Hodge : 

\medskip

\defi Le \textit{groupe de Hodge} d'une vari\'et\'e ab\'elienne $A$, d\'efinie sur un corps de nombres $K$ plong\'e dans $\C$, est not\'e $\Hdg(A)$ et est d\'efini par la formule
\[\Hdg(A)=\MT(A)\cap\SL_V\ \ \text{o\`u }\ \ V=H^1(A(\C),\Q).\]

\medskip

\noindent Les groupes de Mumford-Tate et de Hodge sont reli\'es par $\MT(A)=\G_m\cdot \Hdg(A)$ o\`u $\cdot$ d\'enote le produit presque direct (\textit{i.e.} le morphisme, donn\'e par le produit, du produit direct $\G_m\times \Hdg(A)$ vers $\MT(A)$ est une isog\'enie). 

\medskip

\begin{lemme}\textnormal{\textbf{(Imai \cite{imai})}}\label{produit} Soit $A=\prod_{i=1}^rE_i$ un produit de courbes elliptiques deux \`a deux non-isog\`enes, d\'efinies sur un corps de nombres $K$ plong\'e dans $\C$. On a
\[\Hdg(A)\simeq \prod_{i=1}^r\Hdg(E_i).\]
\end{lemme}

\medskip

\noindent Rappelons maintenant la c\'el\`ebre conjecture de Mumford-Tate. Notons $A/K$ une vari\'et\'e ab\'elienne quelconque de dimension $g\geq 1$, d\'efinie sur un corps de nombres $K$ que l'on suppose plong\'e dans $\C$. Notons \'egalement $\MT(A)$ le groupe de Mumford-Tate correspondant. 

\medskip

\defi\label{rol} Notons $T_{\ell}(A)$ le module de Tate et $V_{\ell}:=T_{\ell}(A)\otimes_{\Z_{\ell}}\Q_{\ell}$. Soit $\ell$ un premier et $\rho_{\ell} : \Gal(\overline{K}/K)\rightarrow \Aut(T_{\ell}(A))\subset\GL(V_{\ell})$ la repr\'esen\-ta\-tion $\ell$-adique associ\'ee \`a l'action de Galois sur les points de $\ell^{\infty}$-torsion de $A$. On d\'efinit $\overline{G_{\ell}}$ (que l'on notera aussi \'eventuellement $G_{\ell,A}$ si besoin est) comme \'etant l'adh\'erence de Zariski de l'image $G_{\ell}$ de $\rho_{\ell}$ dans le groupe alg\'ebrique $\GL_{V_{\ell}}\simeq \GL_{2g,\Q_{\ell}}$. C'est un groupe alg\'ebrique sur $\Q_l$ dont on notera $\overline{G_{\ell}}^0$ la composante neutre (composante connexe de l'identit\'e).

\medskip

\noindent Les th\'eor\`emes de comparaison entre cohomologie \'etale et cohomologie classique (\textit{cf.} \cite{sga4} XI) d'une part, et la comparaison entre le premier groupe de cohomologie \'etale et le module de Tate (\textit{cf.} \cite{milne} 15.1 (a)) d'autre part donnent pour tout premier $\ell$ l'isomorphisme canonique : 
\[V_{\ell}=T_{\ell}(A)\otimes_{\Z_{\ell}}\Q_{\ell}\simeq V\otimes_{\Q}\Q_{\ell}.\]
\noindent Nous fixons une fois pour toute dans la suite un tel isomorphisme. Ceci permet de comparer $\MT(A)\times_{\Q}\Q_{\ell}$ et $\overline{G_{\ell}}^0$. Dans le papier \cite{mumfordtate} (en partie \'ecrit avec Tate comme l'indique Mumford en introduction), se trouve formul\'ee la c\'el\`ebre conjecture~:

\begin{conj}\label{mt}\textnormal{\textbf{(Mumford-Tate)}} Pour tout $\ell$ premier, on a $\overline{G_{\ell}}^0=\MT(A)\times_{\Q}\Q_{\ell}$.
\end{conj}

\medskip

\defi \label{repmt} Notons $\rho : \MT(A)\hookrightarrow \GL_V$ la repr\'esentation naturelle du groupe de Mumford-Tate.

\medskip

\noindent Un certain nombre de cas particuliers, ainsi que de r\'esultats en direction de la conjecture pr\'ec\'edente sont connus. Nous renvoyons \`a la r\'ef\'erence \cite{pink1} pour une discussion d\'etaill\'ee de ces r\'esultats. Disons simplement ici que, suite aux travaux de Serre, cette conjecture est un th\'eor\`eme pour les produits de courbes elliptiques (cf. th\'eor\`eme \ref{MTfort} ci-dessous). De mani\`ere g\'en\'erale on sait par les travaux de Borovo{\u\i} \cite{boro}, Deligne \cite{del} Exp I, 2.9, 2.11, et Pjatecki{\u\i}-{\v{S}}apiro \cite{piat} qu'une inclusion est toujours vraie~:

\medskip

\begin{theo}\label{bdp}\textnormal{\textbf{( Borovo{\u\i}, Deligne, Pjatecki{\u\i}-{\v{S}}apiro)}} Pour tout $\ell$ premier, $\overline{G_{\ell}}^0\subset \MT(A)\times_{\Q}\Q_{\ell}$.
\end{theo}

\medskip 

\noindent Par ailleurs, on a le r\'esultat suivant, d\^u \`a Serre \cite{college8485} 2.2.3. (\textit{cf.} \'egalement \cite{serken}), dans le cas g\'en\'eral.

\medskip

\begin{theo} \textnormal{\textbf{(Serre)}} L'application 
\[\epsilon~: \Gal(\overline{K}/K)\longrightarrow \overline{G_{\ell}}(\Q_{\ell})\longrightarrow\overline{G_{\ell}}(\Q_{\ell})/\overline{G_{\ell}}^0(\Q_{\ell})\]
\noindent est continue surjective, de noyau ind\'ependant de $\ell$ pour tout premier $\ell$.
\end{theo}

\medskip

\noindent Ainsi, le noyau de $\epsilon$ est un sous-groupe d'indice fini de $\Gal(\overline{K}/K)$ donc il existe une extension finie $K'/K$ telle que $\Gal(\overline{K'}/K')\subset \ker\epsilon$. Dit autrement, les deux th\'eor\`emes pr\'ec\'edents donnent le lien entre la repr\'esentation $\rho$ et les repr\'esentations $\ell$-adiques $\rho_{\ell}$ pour tout premier $\ell$ : quitte \`a remplacer au d\'epart $K$ par une extension $K'$ finie ne d\'ependant que de $A$ (ce que nous ferons dans la suite), on a, pour tout premier $\ell$, la factorisation
\[\rho_{\ell} : \Gal(\overline{K}/K)\rightarrow \MT(A)(\Z_{\ell})\overset{\rho}{\rightarrow}\GL_V(\Z_{\ell})\simeq \GL (T_{\ell}(A)).\]

\noindent Ainsi, dans ce contexte la conjecture de Mumford-Tate se traduit en disant que le sous-groupe $G_{\ell}$ est d'indice fini dans $\MT(A)(\Z_{\ell})$. Dans \cite{serrekyoto} conjecture C.3.7.a) (ce qui est not\'e $\Hdg(\Z_{\ell})$ dans cette r\'ef\'erence est ce que nous notons ici $\MT(A)(\Z_{\ell})$), Serre a formul\'e une version encore plus optimiste~:

\begin{conj}\label{mt+}\textnormal{\textbf{(Mumford-Tate, version forte)}} Pour tout $\ell$ premier, on a $G_{\ell}\subset\MT(A)(\Z_{\ell})$, et de plus cette inclusion est d'indice fini, born\'e ind\'ependamment  de $\ell$.
\end{conj}

\noindent Cette version forte de la conjecture de Mumford-Tate est d\'emontr\'ee dans un certain nombre de cas, notamment dans le cas des courbes elliptiques~:

\medskip

\begin{theo}\label{trr}\textnormal{\textbf{(Serre)}} Soit $E/K$ une courbe elliptique. Pour tout $\ell$ premier on a
\[G_{\ell}=\im(\rho_{\ell})\subset \MT(E)(\Z_{\ell}),\]
\noindent cette inclusion \'etant de conoyau fini, born\'e ind\'ependamment de $\ell$. Si $E$ est sans multiplication complexe, l'inclusion est m\^eme une \'egalit\'e pour tout $\ell$ assez grand (d\'ependant de $E$).
\end{theo}
\demo Dans le cas des courbes elliptiques de type CM cela d\'ecoule du th\'eor\`eme 1 p.II-26 de \cite{serreabelien}. Dans le cas des courbes elliptiques sans multiplication complexe il s'agit du th\'eor\`eme 3 p. 299 de \cite{serregalois}.\hfill$\Box$

\medskip

\noindent Mieux, Serre donne dans \cite{serregalois} une preuve de cette conjecture forte pour un produit de deux courbes elliptiques non-isog\`enes, Ribet  \cite{ribetmodulaire} l'en d\'eduit pour un produit de courbes elliptiques deux \`a deux non-isog\`enes sans multiplication complexe, et, dans une lettre \`a Masser (cf. \cite{serremasser}), Serre donne une preuve pour un produit de courbes elliptiques deux \`a deux non-isog\`enes \`a multiplication complexe. Il est alors essentiellement formel d'en d\'eduire le r\'esultat pour un produit de courbes elliptiques~:

\medskip

\begin{theo}\label{MTfort} Soient $g\geq 1$ et $E_1,\ldots, E_g/K$ des courbes elliptiques. Posons $A=\prod_{i=1}^gE_i$. Pour tout $\ell$ premier on a
\[G_{\ell,A}=\im(\rho_{\ell})\subset \MT(A)(\Z_{\ell}),\]
\noindent cette inclusion \'etant de conoyau fini, born\'e ind\'ependamment de $\ell$. 
\end{theo}
\demo Notons tout d'abord qu'il suffit de prouver le r\'esultat dans le cas o\`u les $E_i$ sont deux \`a deux non-isog\`enes, le cas g\'en\'eral d'un produit $\prod_{i=1}^gE_i^{n_i}$ en d\'ecoule, car par le lemme \ref{puissance0} le groupe $\MT(\prod_{i=1}^gE_i^{n_i})=\MT(\prod_{i=1}^g E_i)$ via le plongement diagonal des $\GL_{V_i}$ dans $\GL_{V_i^{\oplus n_i}}$ (notations du lemme \ref{puissance0}) et il en est de m\^eme pour l'image  des repr\'esentations $\ell$-adiques $\rho_{\ell}$. Si toutes les $E_i$ sont \`a multiplications complexes, le r\'esultat est prouv\'e par Serre dans \cite{serremasser}. Si au contraire aucune $E_i$ n'a de multiplication complexe, le r\'esultat est le th\'eor\`eme 6 de Serre \cite{serregalois} si $g=2$, ce r\'esultat \'etant \'etendu \`a $g$ quelconque par Ribet \cite{ribetmodulaire} Theorem 3.5. Il reste donc \`a traiter le cas o\`u au moins l'une des $E_i$ est \`a multiplications complexes, sans que toutes le soit. L\`a encore le cas $g=2$ est prouv\'e dans \cite{serregalois} : il s'agit du th\'eor\`eme 7. La preuve donn\'ee dans \cite{serregalois} pour $g=2$ reli\'ee \`a l'argument de Ribet permettent l\`a encore d'\'etendre le r\'esultat au cas $g$ quelconque : \'Ecrivons le produit $\prod_{i=1}^gE_i$ sous la forme $B\times C$ o\`u $B=\prod_{i=1}^rB_i$ est un produit non-vide de courbes elliptiques sans multiplications complexes, et o\`u $C$ est un produit non-vide de courbes elliptiques ayant des multiplications complexes. Par la conjecture de Mumford-Tate forte pour $B$ et $C$, valable par ce qui pr\'ec\`ede, on sait que l'on a
\[\Hdg(B)(\Z_{\ell})\cong \Gal\left(K(B[\ell^{\infty}])/K(\mu_{\ell^{\infty}})\right)\ \ \text{ et de m\^eme pour $C$,}\]
\noindent o\`u $\cong$ signifie, \`a indice fini born\'e ind\'ependamment de $\ell$ pr\`es. Le groupe de Hodge d'un produit de courbes elliptiques \'etant isomorphe au produit des groupes de Hodge, on a ainsi
\[\Hdg(A)(\Z_{\ell})\cong\Gal\left(K(B[\ell^{\infty}])/K(\mu_{\ell^{\infty}})\right)\times\Gal\left(K(C[\ell^{\infty}])/K(\mu_{\ell^{\infty}})\right).\]
\noindent Il nous suffit donc de prouver que les extensions $K(B[\ell^{\infty}])$ et $K(C[\ell^{\infty}])$ sont presque lin\'eairement disjointes au dessus de $K(\mu_{\ell^{\infty}})$, \textit{i.e.}, montrer que l'extension $K_{B,C}:=K(B[\ell^{\infty}])\cap K(B[\ell^{\infty}])$ est finie sur $K(\mu_{\ell^{\infty}})$, de degr\'e born\'e ind\'ependamment de $\ell$. Nous suivons pour cela la preuve du th\'eor\`eme 6 de \cite{serregalois} : l'extension $K_{B,C}$ est ab\'elienne sur $K(\mu_{\ell^{\infty}})$ car $C$ est une vari\'et\'e ab\'elienne de type CM. Notons $\Kc$ l'extension de $K$ engendr\'ee par les $K(\mu_{\ell^{\infty}})$ lorsque $\ell$ varie et, pour toute vari\'et\'e ab\'elienne $X/K$, notons $K(X_{\infty})$ l'extension de $K$ engendr\'ee par tous les points de torsion de $X(\overline{K})$.

\noindent Nous voulons montrer que  $\Gal(K(B_{\infty})/\Kc)$ contient $\Gal(K(B_{\infty})/K_{B,C})$ comme sous-groupe ouvert. Or a un indice fini pr\`es, on sait que \[\Gal(K(B_{\infty})/\Kc)=\prod_{\ell}\prod_{i=1}^r\Hdg(B_i)(\Z_{\ell})\cong\prod_{\ell}\prod_{i=1}^r\SL_2(\Z_{\ell}).\]
\noindent Comme indiqu\'e en remarque 3. du lemme 3.4 de \cite{ribetmodulaire}, chacun des groupes $\prod_{\ell}\SL_2(\Z_{\ell})$ satisfait la propri\'et\'e de commutateurs exig\'ee pour l'application du lemme 3.4 de Ribet \cite{ribetmodulaire} : pour tout sous groupe ouvert $U$ de $\prod_{\ell}\SL_2(\Z_{\ell})$, l'adh\'erence du sous-groupe des commutateurs de $U$ est ouverte dans $\prod_{\ell}\SL_2(\Z_{\ell})$. Pour pouvoir appliquer le lemme 3.4 de \cite{ribetmodulaire} (et ainsi conclure la preuve), nous reste donc \`a v\'erifier que la projection de $\rho_{\ell}\left(\Gal(K(B_{\infty})/K_{B,C})\right)$ sur chacun des double facteur $\prod_{\ell}(\SL_2\times \SL_2)(\Z_{\ell})$ est ouverte. L\`a encore on suit l'argument de Serre : $\SL_2\times\SL_2$ admet $\ssl_2\times \ssl_2$ comme alg\`ebre de Lie. Cette derni\`ere est \'egale \`a son alg\`ebre de Lie d\'eriv\'ee, donc la remarque 2 suivant le lemme 3.4 de \cite{ribetmodulaire} nous assure que $\prod_{\ell}(\SL_2\times\SL_2)(\Z_{\ell})$ v\'erifie la propri\'et\'e de commutateurs rappel\'ee ci-dessus. On l'applique alors avec comme groupe $U$, l'image de $\Gal(K(B_{\infty})/\Kc)$ dans $\prod_{\ell}(\SL_2\times\SL_2)(\Z_{\ell})$, image qui est ouverte par le th\'eor\`eme 6' de \cite{serregalois}. Ceci permet d'appliquer le lemme 3.4 de \cite{ribetmodulaire} et donc de conclure.\hfill$\Box$

\medskip

\noindent Nous rappelons enfin, un fait indiqu\'e dans l'introduction, sans doute bien connu des experts, mais que nous n'avons pas trouv\'e dans la litt\'erature.

\medskip

\begin{prop}\label{mtmin} Soit $A$ une vari\'et\'e ab\'elienne sur $\overline{\Q}$. Le groupe de Mumford-Tate $\MT(A)$ est de dimension $2$ si et seulement si $A$ est isog\`ene \`a une puissance d'une courbe elliptique \`a multiplication complexe.
\end{prop}

\demo Dire que le groupe de Mumford-Tate de $A$ est de dimension $2$ \'equivaut \`a dire que son groupe de Hodge, $\Hdg(A)$ est de dimension $1$. Or tout groupe alg\'ebrique connexe de dimension $1$ sur $\Q$ est commutatif. Ceci entra\^{\i}ne (\textit{cf.} par exemple \cite{mumford} paragraphe 2, seconde proposition) que $A$ est une vari\'et\'e ab\'elienne de type CM (il s'agit m\^eme d'une propri\'et\'e \'equivalente). Or si $A$ est simple de type CM, on peut utiliser la borne de Ribet (\cite{ribet}) : 
\[\dim \Hdg(A)\geq 1+\log_2(\dim A).\]
\noindent Donc si $A$ est simple, $\dim\Hdg(A)=1$ implique que $A$ est une courbe elliptique de type CM. Si $A$ est isog\`ene \`a un produit $\prod_{i=1}^n A_i^{n_i}$ avec les $A_i$ deux \`a deux non-isog\`enes, alors pour tout $i$, 
\[\Hdg(A)=\Hdg(A_1\times\ldots\times A_n)\twoheadrightarrow \Hdg(A_i).\]
\noindent Ceci entra\^ine que si $\dim\Hdg(A)=1$ alors les $A_i$ sont de dimension $1$, donc, en appliquant par exemple le lemme \ref{produit} pr\'ec\'edent, on voit que 
 $\Hdg(A)=\prod_{i=1}^n\Hdg(A_i)$ et donc que n\'ecessaire\-ment $n=1$. Finalement, $A=A_1^{n_1}$ avec $A_1$  simple de type CM et de dimension $1$. Ainsi $A$ est bien  une puissance d'une courbe elliptique \`a multiplication complexe. La r\'eciproque est claire.\hfill$\Box$

\medskip

\noindent Concluons ce paragraphe en indiquant des hypoth\`eses  garantissant que le groupe de Hodge d'un produit se d\'ecompose en le produit des groupes de Hodge. Notons que l'on a toujours
$\Hdg(A_1\times\dots\times A_n)\subset\Hdg(A_1)\times\dots\times \Hdg(A_n)$ mais que l'\'egalit\'e n'est pas toujours vraie ({\it cf}. par exemple \cite{mz} paragraphes 3 et 5.3). Rappelons pour cela la notion de dimension relative~:  

\medskip

\defi Soit $A/\C$ une vari\'et\'e ab\'elienne simple. Notons $D=\End(A)\otimes\Q$ et $E$ le centre de $D$. La \textit{dimension relative de $A$}, not\'ee $\dim_{\text{rel}}(A)$, est par d\'efinition l'entier donn\'e par

\[\dim_{\text{rel}}(A):=\begin{cases}
\frac{\dim A}{[E:\Q]}\			& 					\text{ si $A$ est de type I au sens de la classification d'Albert}\\ 
\frac{\dim A}{2[E:\Q]}\ 		&					 \text{ si $A$ est de type II ou III}\\
\frac{2\dim A}{\sqrt{[D:E]}[E:\Q]}\ &	\text{ si $A$ est de type IV.}
\end{cases}\]

\medskip

\begin{theo}\textnormal{\textbf{(Ichikawa \cite{ichi})}} \label{prodhdg}Soit $A/\mathbb{C}$ une vari\'et\'e ab\'elienne dont les sous-vari\'et\'es ab\'eliennes simples sont toutes de dimension relative impaire. Soient $B$ et $C$ deux vari\'et\'es ab\'eliennes telles que $A$ est isog\`ene \`a $B\times C$ et telles que $B$ (respectivement $C$) est de type I, II ou III (respectivement IV) au sens de la classification d'Albert. Alors
\[\Hdg(A)=\Hdg(B)\times\Hdg(C).\]
\noindent De plus, avec les notations pr\'ec\'edentes, si $B$ est isog\`ene au produit $\prod_{i=1}^rB_i^{n_i}$, les $B_i$ \'etant deux \`a deux non-isog\`enes, alors
\[\Hdg(B)=\prod_{i=1}^r\Hdg(B_i).\]
\end{theo}

\section{Preuve de la proposition \ref{p1}\label{pp1}}

\noindent Soient $I$ un ensemble non-vide de $\{1,\ldots,n\}$ et $F/K_0$ une extension finie. Notons $A_I:=\prod_{i\in I}A_i^{n_i}$. Pour tout $\epsilon>0$, on a
\[\left|\left(A_I\right)(F)_{\tors}\right|\leq \left|\left(\prod_{i=1}^nA_i^{n_i}\right)(F)_{\tors}\right|=\left|A(F)_{\tors}\right|\ll[F:K_0]^{\gamma(A)+\epsilon}.\]

\medskip

\noindent On en d\'eduit que, pour tout tel ensemble $I$, on a $\gamma(A)\geq\gamma(A_I)$. Par ailleurs le th\'eor\`eme 1.4 de \cite{ratazziens} assure que 
\[\gamma(A_I)\geq \frac{2 \dim A_I}{\dim \MT(A_I)}=\frac{2\sum_{i\in I}n_i\dim A_i}{\dim \MT(A_I)}.\]
\noindent Or sur la d\'efinition du groupe de Mumford-Tate (voir le lemme \ref{puissance0}), on voit que 
\[\MT(A_I)=\MT\left(\prod_{i\in I}A_i^{n_i}\right)\simeq\MT\left(\prod_{i\in I}A_i\right).\]
\noindent On d\'eduit donc l'in\'egalit\'e voulue.\hfill$\Box$

\medskip

\rem Le m\^eme argument (loc. cit.) indique en fait que, pour toute vari\'et\'e ab\'elienne $A/K$, on a en fait $\gamma(A)\geq \frac{2\dim A}{d_{\ell}}$, o\`u $d_{\ell}:=\dim \overline{G_{\ell}}$. On voit ainsi que si la conjecture de Mumford-Tate est fausse pour $A$ et $\ell$, la r\'eponse \`a la question \ref{conjR}  sera n\'egative.

\medskip

\noindent Notons que ce r\'esultat est d\'ej\`a suffisant pour montrer que la question na\"ive que l'on aurait pu poser, \`a savoir  ``a-t-on $\gamma(A)=\frac{2\dim A}{\dim \MT(A)}$ ?" admet une r\'eponse n\'egative en g\'en\'eral. En effet, il suffit pour cela de trouver un exemple de vari\'et\'e ab\'elienne $A$ telle que $\alpha(A)$ est strictement plus grand que $\frac{2\dim A}{\dim \MT(A)}$. La seconde partie de la proposition suivante donne le contre-exemple recherch\'e.

\medskip

\begin{prop}Soient $n\geq 1$ un entier et $E_1,\ldots,E_n/\overline{\Q}$ des courbes elliptiques deux \`a deux non-isog\`enes. Posons $A=\prod_{i=1}^nE_i$. Si toutes les $E_i$ sont \`a multiplication complexe, ou si aucune $E_i$ est \`a multiplication complexe, alors 
\[\alpha(A)=\frac{2\dim A}{\dim \MT(A)}=\begin{cases}
\frac{2n}{1+3n}	&	\text{ dans le cas sans multiplication complexe}\\
\frac{2n}{1+n}	&	\text{ dans le cas avec multiplication complexe}.\\
\end{cases}\]
\noindent Sinon $A=\prod_{i=1}^rE_i\times\prod_{i=r+1}^nE_i$ avec $n>r\geq 1$, les $E_i$, $1\leq i\leq r$ \`a multiplication complexe et les $E_i$, $r+1\leq i\leq n$ sans multiplication complexe. En posant $B=\prod_{i=1}^rE_i$, on a
\[\alpha(A)=\frac{2\dim B}{\dim\MT(B)}=\frac{2r}{1+r}>\frac{2n}{1+r+3(n-r)}=\frac{2\dim A}{\dim \MT(A)}.\]
\end{prop}
\demo C'est un simple calcul, qui repose sur les faits suivants~: 
\begin{enumerate}
\item Si $E$ est une courbe elliptique \`a multiplication complexe, alors $\dim\MT(E)=2$ ; 
\item si $E$ est une courbe elliptique sans multiplication complexe, alors $\dim \MT(E)=4$ ; 
\item si $X=\prod_{i=1}^m E_i$ est un produit de courbes elliptiques deux \`a deux non-isog\`enes, alors $\Hdg(X)\simeq \prod_{i=1}^m\Hdg(E_i)$ d'apr\`es le lemme \ref{produit}.
\end{enumerate}
\noindent Ces trois points permettent ais\'ement de conclure. \hfill$\Box$

\section{Quelques r\'eductions}

\noindent Nous commen\c{c}ons par deux r\'eductions (au cas $\ell$-adique et au cas d'un groupe produit) valables en toute g\'en\'eralit\'e pour des vari\'et\'es ab\'eliennes quelconques.

\subsection{R\'eduction au cas $\ell$-adique}

\noindent Soit $A/\overline{\Q}$ une vari\'et\'e ab\'elienne d\'efinie sur un corps de nombres $K$. 

\begin{prop}\label{ladique}Soit $\alpha>0$. Pour d\'emontrer que $\gamma(A)\leq\alpha$, il suffit de montrer que : il existe une cons\-tan\-te strictement positive $C(A/K)$ ne d\'ependant que de $A/K$ telle que pour tout nombre premier $\ell$, pour tout sous-groupe fini $H_{\ell}$ de $A[\ell^{\infty}]$, on a
\begin{equation}\label{equa}
\text{Card}\left(H_{\ell}\right)\leq C(A/K)[K(H_{\ell}):K]^{\alpha}.
\end{equation}
\end{prop}
\demo Soit $L/K$ une extension finie. Posons
\[ H=A(L)_{\tors}\ \text{ et pour tout premier $\ell$, }\ H_{\ell}=A(L)_{\tors}\left[\ell^{\infty}\right].\]
\noindent Nous pouvons appliquer l'hypoth\`ese de l'\'enonc\'e aux groupes $H_{\ell}$ pour tout premier $\ell$~:
\[\left| H_{\ell} \right|\ll [K(H_{\ell}):K]^{\alpha}.\]
\noindent En notant $\omega(n)$ le nombre de nombre premiers divisant $n$, il vient 
\begin{equation}\label{in}
\left|A(L)_{\tors}\right|=|H|=\prod_{\ell}|H_{\ell}|\leq C(A/K)^{\omega\left(|A(L)_{\tors}|\right)}\prod_{\ell}[K(H_{\ell}):K]^{\alpha}.
\end{equation}
\noindent Par ailleurs, une estimation classique de $\omega(n)$ est la suivante (\textit{cf.} par exemple \cite{tenenbaum} p. 85 \S\ 5.3)~: $\omega(n)\ll\frac{\log n}{\log\log n}$. En se souvenant ({\it confer} par exemple le th\'eor\`eme \ref{theomasser}) que $\left|A(L)_{\tors}\right|\ll[L:K]^c$, on en tire
\[\omega(\left|A(L)_{\tors}\right|)\ll\frac{\log\left|A(L)_{\tors}\right|}{\log\log\left|A(L)_{\tors}\right|}\ll\frac{\log[L:K]}{\log\log[L:K]}\]
 L'in\'egalit\'e (\ref{in}) peut donc se r\'e\'ecrire
\[\left|A(L)_{\tors}\right| \leq C^{\frac{\log[L:K]}{\log\log[L:K]}}\prod_{\ell}[K(H_{\ell}):K]^{\alpha}.\]
\noindent De plus un th\'eor\`eme de Serre (Th\'eor\`eme 1 de \cite{college8586}) assure que les repr\'esentations $\rho_{\ell} : G_K\rightarrow \Aut(T_{\ell}(A))$ sont ind\'ependantes (au moins quitte \`a \^etre mont\'e au d\'epart, sur une extension finie de $K$ ne d\'ependant que de $A/K$, ce que l'on suppose ici). Comme pour tout premier $\ell$, $H_{\ell}\subset A[\ell^{\infty}]$, on en d\'eduit que
 $\prod_{\ell}[K(H_{\ell}):K]\ll[K(H):K]\leq[L:K]$ et enfin que, pour tout $\epsilon>0$, on a bien~:
\[\left|A(L)_{\tors}\right|\leq C^{\frac{\log[L:K]}{\log\log[L:K]}}[K(H):K]^{\alpha}\leq [L:K]^{\alpha+\epsilon}.\]
\noindent Ceci conclut.\hfill$\Box$

\subsection{R\'eduction au cas produit}

\noindent Soit $A=\prod_{i=1}^rA_i$ un produit de vari\'et\'es ab\'eliennes quelconques (isog\`enes ou non). Soient $\ell$ un nombre premier et $H$ un sous-groupe fini de $A[\ell^{\infty}]$. Notons $H'=\prod_{i=1}^r\pr_i(H)$ o\`u $\pr_i : A\rightarrow A_i$ est la projection canonique. Si $(P_1,\ldots,P_r)\in H'$ alors chaque $P_i$ est d\'efini sur $K(H)$ (car pour tout $i$ entre $1$ et $r$, il existe $Q_j$, $j\not=i$ tels que $(Q_1,\ldots,Q_{i-1},P_i,Q_{i+1},\ldots,Q_r)\in H$). Or $H\subset H'$, donc $K(H)=K(H')$. Ceci permet de se ramener au cas d'un groupe produit : si $H'$ v\'erifie $|H'|\ll [K(H'):K]^{\star}$, il en est de m\^eme, avec le m\^eme exposant $\star$, pour $H$. Ainsi dans toute la suite nous pourrons supposer que $H$ est un sous-groupe produit. Mieux : 

\noindent Si $A=\prod_{i=1}^nA_i^{n_i}$ les $A_i$ \'etant deux \`a deux non-isog\`enes sur $\overline{\Q}$, si $\ell$ est un premier quelconque et si $H=\prod_{i=1}^n(\prod_{j=1}^{n_i}H_{i,j})$ o\`u les $H_{i,j}$ sont des sous-groupes finis de $A_i[\ell^{\infty}]$ pour tout $i$ et $j$~; posons $L=K(H)$, on a alors $H\subset \prod_{i=1}^nA_i(L)_{\rm tors}^{n_i}$. De plus l'extension engendr\'ee sur $K$ par les $A_i(L)_{\tors}$ est incluse dans $L=K(H)$.

\noindent Ceci montre que, si $A$ est de la forme $\prod_{i=1}^nA_i^{n_i}$,  on peut supposer dans la suite le groupe $H$ de $A[\ell^{\infty}]$ de la forme $\prod_{i=1}^nH_i^{n_i}$, o\`u $H_i$ est un sous-groupe fini de $A_i[\ell^{\infty}]$.

\subsection{Une r\'eduction sp\'ecifique au cas de type CM\label{redgalois}}
\noindent Supposons dans ce paragraphe que la vari\'et\'e ab\'elienne $A$ est de type CM. Soit $\ell$ un premier et $\H\subset A[\ell^{\infty}]$ fini. L'extension $K(H)/K$ est ab\'elienne donc galoisienne, de groupe de Galois $\mathcal{G}_{\ell}$. Si $x\in H$ et $\sigma\in\mathcal{G}_{\ell}$, alors $\sigma(x)$ est d\'efini sur $K(H)$. Donc le groupe $H_{\mathcal{G}_{\ell}}:=\mathcal{G}_{\ell}\cdot H$ engendr\'e par les $\sigma(x)$ est d\'efini sur $K(H)$ et contient $H$, donc $K(H)=K(H_{\mathcal{G}_{\ell}})$. Ainsi l'in\'egalit\'e 
\[\left|H_{\mathcal{G}_{\ell}}\right|\ll\left[K\left(H_{\mathcal{G}_{\ell}}\right):K\right]^{\star}\]
\noindent entra\^ine l'in\'egalit\'e
\[ |H|\ll [K(H):K]^{\star}.\]
\noindent Ainsi dans le cas CM nous pourrons toujours supposer que le groupe $H$ est Galoisien.

\section{Pr\'eliminaires  galoisiens}

\noindent R\'esumons ici les arguments galoisiens que nous allons utiliser de mani\`ere r\'ep\'et\'ee.

\begin{lemme}\label{galelem}
\noindent Soit $K/k$ une extension galoisienne, $F/k$ une extension quelconque, alors $KF/F$ est galoisienne et la restriction au corps $K$ induit un isomorphisme de $\Gal(KF/F)$ sur $\Gal(K/F\cap K)$.
\end{lemme} 
\demo
Voir Lang, {\it Algebra},\cite{algebra} [Ch. VIII, \S 1, Theorem 4] 
$\findemo$

\medskip

\noindent Remarquons que, si $k=K\cap F$ on obtient $[KF:K]=[F:k]$ et, par cons\'equent, $[K:F]=[K:k]$.

\medskip

\defi Nous dirons que les corps $k$, $K_1$, $K_2$ et $L$ forment un {\it parall\'elogramme} si $L$ est le compositum de $K_1$ et $K_2$ et si $[L:k]=[K_1:k][K_2:k]$ (cette derni\`ere condition \'equivaut \`a
$[L:K_1]=[K_2:k]$ ou encore $[L:K_2]=[K_1:k]$.
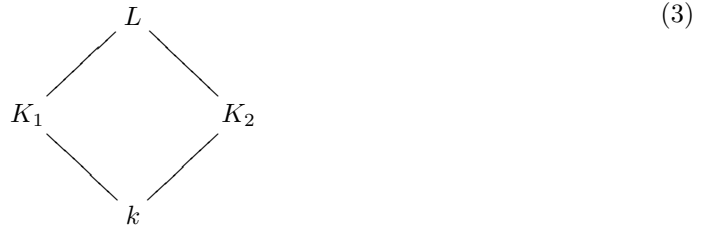
 \begin{equation}\label{exparal}
\xymatrix{
&L\ar@{-}[dl] \ar@{-}[dr]&\\
 K_1\ar@{-}[dr]& & K_2\ar@{-}[dl] \\
 &k&}
 \end{equation}
 On voit ais\'ement que si $L=K_1K_2$, une condition n\'ecessaire pour qu'on ait un parall\'elogramme est que $k=K_1\cap K_2$. Le lemme \ref{galelem} indique que, si l'une des deux extensions $K_1/k$ ou $K_2/k$ est galoisienne cette condition est aussi suffisante. On peut compl\'eter cette discussion avec le lemme suivant.
 
 \begin{lemme} Supposons que $k$, $L_1$, $L_2$ et $L=L_1L_2$ forment un parall\'elogramme. Soit
 $k\subset K_1\subset L_1$ et $k\subset K_2\subset L_2$, alors $k$, $K_1$, $K_2$ et $K=K_1K_2$ forment un parall\'elogramme.  
 
 \end{lemme}
 \demo La situation est r\'esum\'ee dans le diagramme suivant~:
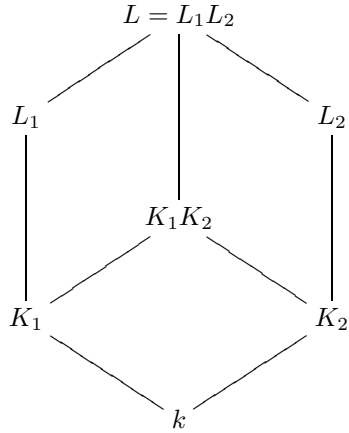
  $$
\xymatrix{
& L=L_1L_2\ar@{-}[dl]  \ar@{-}[dd] \ar@{-}[dr] & \\
L_1\ar@{-}[dd]&&L_2\ar@{-}[dd] \\
&K_1K_2\ar@{-}[dl] \ar@{-}[dr]&\\
 K_1\ar@{-}[dr]& & K_2\ar@{-}[dl] \\
 &k&}
 $$
 Le fait que $[L:L_2]=[L_1:k]$ (et sym\'etriquement $[L:L_1]=[L_2:k]$) permet d'\'ecrire
 $$[L_1L_2:K_1K_2]=[L_1L_2:L_1K_2][L_1K_2:K_1K_2]\leq [L_2:K_2][L_1:K_1].$$
 Si l'on pose $c:=[K_1:k][K_2:k]/[K_1K_2:k]$ (qui est $\geq 1)$) on a~:
 $$[L_1L_2:k]= [L_1L_2:K_1K_2][K_1K_2:k]\leq[L_2:K_2][L_1:K_1][K_1:k][K_2:k]/c=[L_1:k][L_2:k]/c$$
 On en tire donc $c=1$ et la conclusion voulue.
 $\findemo$

\medskip

\noindent Remarquons que ce lemme permet de d\'efinir les parall\'elogrammes infinis : les corps $k$, $L_1$, $L_2$ et  le compositum $L=L_1L_2$ forment un parall\'elogramme (avec $L_i/k$ \'eventuellement infinie) si pour toutes $K_i$, sous-extensions finies sur $k$, les extensions  $k$, $K_1$, $K_2$ et $K=K_1K_2$ forment un parall\'elogramme. 

\section{Sous-extensions  de $K\left(A\left[\ell^{\infty}\right]\right)/K$\label{conv}}

\subsection{Produit de vari\'et\'es ab\'eliennes}

\noindent {\bf Convention.} {\it Nous utiliserons dans ce paragraphe la convention suivante~: on \'ecrira $[L:K]=[F:K]$ au lieu de $c_1^{-1}[F:k]\leq [L:K]\leq c_1[F:k]$ et de m\^eme on s'autorisera \`a \'ecrire $L=K$ si on a $\max\{[L:L\cap K],[K:L\cap K]\}\leq c_1$. (La constante $c_1$ ne d\'ependant \'evidemment que de la vari\'et\'e ab\'elienne $A$ et du corps de base $K$.)} 

\medskip

\noindent Soient $A_1$ et $A_2$ deux vari\'et\'es ab\'eliennes d\'efinies sur $K$.
Consid\'erons $A=A_1\times A_2$ et faisons les deux hypoth\`eses suivantes : 
\begin{enumerate}
\item $\hdg(A)\simeq\hdg{A_1}\times\hdg{A_2}$.
\item $A=A_1\times A_2$ v\'erifie la conjecture de Mumford-Tate forte (conjecture \ref{mt+}).
\end{enumerate}

\medskip

\rem Notons que les deux conditions pr\'ec\'edentes, bien que plac\'ees sur le m\^eme plan, n'ont en fait pas du tout le m\^eme statut : on s'attend \`a ce que la condition 2. soit toujours v\'erifi\'ee. Il n'en est pas de m\^eme pour la condition $1.$ dont on sait au contraire qu'elle peut \^etre mise en d\'efaut (comme indiqu\'e au paragraphe \ref{mtt}, \textit{cf.} par exemple \cite{mz} paragraphes 3 et 5.3). Nous avons indiqu\'e dans la proposition \ref{prodhdg} du paragraphe \ref{mtt} certains r\'esultats connus concernant cette condition $1.$

\medskip

\begin{prop}\label{h4} Si $A=A_1\times A_2$ v\'erifie les conditions 1. et 2. pr\'ec\'edentes, alors
\begin{equation}\label{hodgehyp}
\Gal(K(A[\ell^{\infty}])/K(\mu_{\ell^{\infty}}))\cong \Gal(K(A_1[\ell^{\infty}])/K(\mu_{\ell^{\infty}}))\times\Gal(K(A_2[\ell^{\infty}])/K(\mu_{\ell^{\infty}}))
\end{equation}
\noindent (o\`u $\cong$ signifie isomorphe \`a indice fini pr\`es).
\end{prop}
\demo Imm\'ediat.\hfill$\Box$

\medskip

\noindent Consid\'erons maintenant un sous-groupe $H=H_1\times H_2\subset A[\ell^{\infty}]$, avec $H_1$ et $H_2$ finis. On dispose du diagramme suivant d'extensions.
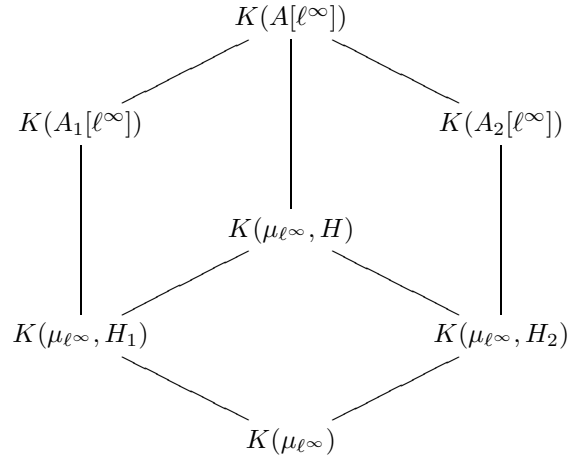
\[
\xymatrix{
& K(A[\ell^{\infty}])\ar@{-}[dl]  \ar@{-}[dd] \ar@{-}[dr] & \\
K(A_1[\ell^{\infty}])\ar@{-}[dd]&&K(A_2[\ell^{\infty}])\ar@{-}[dd] \\
&K(\mu_{\ell^{\infty}},H)\ar@{-}[dl] \ar@{-}[dr]&\\
 K(\mu_{\ell^{\infty}},H_1)\ar@{-}[dr]& & K(\mu_{\ell^{\infty}},H_2)\ar@{-}[dl] \\
 &K(\mu_{\ell^{\infty}})&}
 \]
\noindent On peut utiliser le m\^eme diagramme en un cran fini en rempla\c{c}ant $\ell^{\infty}$ par $\ell^N$ avec $N$ suffisamment grand pour que $H\subset A[\ell^N]$. La proposition \ref{h4} peut se traduire en disant que le diagramme ci-dessus fournit un parall\'elogramme infini.

\medskip

\defi Soit $A$ une vari\'et\'e ab\'elienne sur un corps de nombres $K$. Nous dirons que \textit{$A$ v\'erifie la propri\'et\'e $\mu$} si, pour tout premier $\ell$ et pour tout sous-groupe fini $H$ de $A[\ell^{\infty}]$, il existe un entier $m$ tel que, \`a indice fini born\'e (ind\'ependamment de $\ell$) pr\`es~:
\begin{equation}\label{intercyc}
K(H)\cap K(\mu_{\ell^{\infty}})=K(\mu_{\ell^{m}}).
\end{equation}

\rem Pour les calculs combinatoires ult\'erieurs n\'ecessaires dans le calcul de $\gamma(A)$ pour $A=\prod_{i=1}^rA_i^{n_i}$, il est de plus important de savoir d\'eterminer l'exposant $m$ associ\'e \`a un groupe fini $H$ donn\'e.

\medskip

\begin{prop}\label{casr=2} Soit $A=A_1\times A_2$ v\'erifiant la conclusion de la proposition \ref{h4}. Soit $H=H_1\times H_2$ comme pr\'ec\'edemment. Supposons de plus que $A_1$ et $A_2$ v\'erifient la propri\'et\'e $\mu$ avec un exposant $m_i$ associ\'e au sous-groupe $H_i$ pour $i\in\{1,2\}$. Posons $K_1=K(H_1)\cap K(H_2)$, on a
\begin{equation}
[K(H):K_1]=[K(H_1):K_1][K(H_2):K_1].
\end{equation}
Autrement dit $K_1$, $K(H_1)$, $K(H_2)$ et $K(H)$ forment un parall\'elogramme.
De plus, si on note $m:=\max(m_1,m_2)$, on a l'\'egalit\'e~:
\begin{equation}
K(H)\cap K(\mu_{\ell^{\infty}})=K(\mu_{\ell^m}),
\end{equation}
\noindent autrement dit $A$ v\'erifie la propri\'et\'e $\mu$ avec l'exposant $m$ associ\'e \`a $H$.
\end{prop}
\demo Comme nous l'avons dit, on peut traduire la conclusion de la proposition \ref{h4} (identit\'e (\ref{hodgehyp})) en disant que le diagramme ci-dessous est un parall\'elogramme (en prenant $N$ assez grand pour que $H\subset A[\ell^N]$).
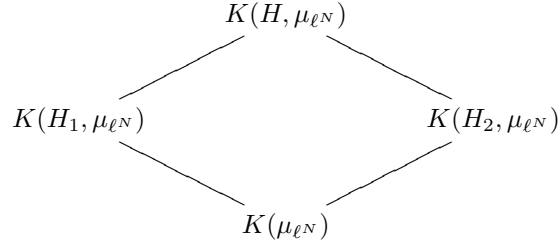
 $$
 \xymatrix{
&K(H,\mu_{\ell^N})\ar@{-}[dl] \ar@{-}[dr]&\\
K(H_1,\mu_{\ell^N})\ar@{-}[dr]& & K(H_2,\mu_{\ell^N})\ar@{-}[dl] \\
 &K(\mu_{\ell^N})&}
 $$
 Supposons (sans perte de g\'en\'eralit\'e) que $m_1\leq m_2$, alors $K_1=K(\mu_{\ell^{m_1}})$; en effet $K(\mu_{\ell^{m_1}})$ est contenu dans $K(H_1)\cap K(H_2)$ qui doit \^etre contenu dans $K(\mu_{\ell^{\infty}})$ et donc dans $K(H_1)\cap K(\mu_{\ell^{\infty}})=K(\mu_{\ell^{m_1}})$.  De m\^eme $K(H_1,\mu_{\ell^{m_2}})\cap K(\mu_{\ell^{N}})=K(\mu_{\ell^{m_2}})$. En effet la propri\'et\'e $\mu$ appliqu\'ee aux $H_i$ se traduit par le fait que le diagramme ci-dessous est un parall\'elogramme.
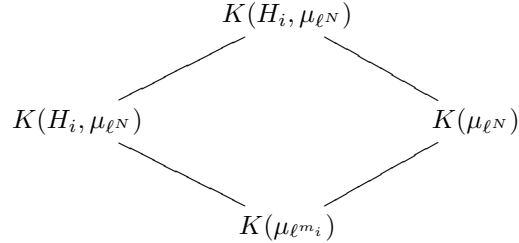
 $$
 \xymatrix{
&K(H_i,\mu_{\ell^N})\ar@{-}[dl] \ar@{-}[dr]&\\
K(H_i,\mu_{\ell^N})\ar@{-}[dr]& & K(\mu_{\ell^N})\ar@{-}[dl] \\
 &K(\mu_{\ell^{m_i}})&}
 $$
On en tire que  le diagramme ci-dessous est un parall\'elogramme (car c'est un sous-parall\'elogramme du diagramme pr\'ec\'edent pour $i=1$)~:
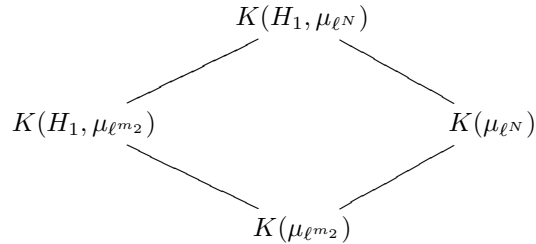
 $$
 \xymatrix{
&K(H_1,\mu_{\ell^N})\ar@{-}[dl] \ar@{-}[dr]&\\
K(H_1,\mu_{\ell^{m_2}})\ar@{-}[dr]& & K(\mu_{\ell^N})\ar@{-}[dl] \\
 &K(\mu_{\ell^{m_2}})&}
 $$
 et en particulier $K(H_1,\mu_{\ell^{m_2}})\cap K(\mu_{\ell^{N}})=K(\mu_{\ell^{m_2}})$. 
 
\noindent Par cons\'equent  le diagramme ci-dessous est aussi un parall\'elogramme.
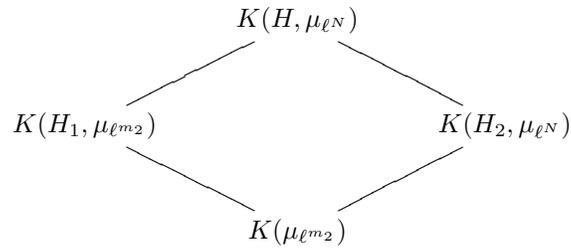
 $$
 \xymatrix{
&K(H,\mu_{\ell^N})\ar@{-}[dl] \ar@{-}[dr]&\\
K(H_1,\mu_{\ell^{m_2}})\ar@{-}[dr]& & K(H_2,\mu_{\ell^{N}})\ar@{-}[dl] \\
 &K(\mu_{\ell^{m_2}})&}
 $$
On obtient ainsi un sous-parall\'elogramme~:
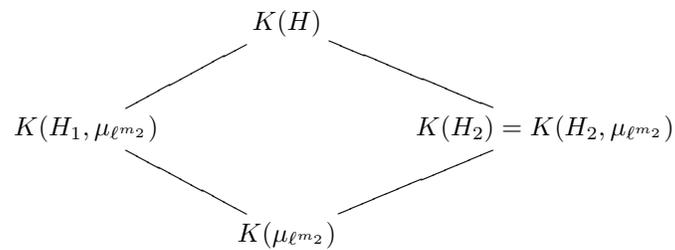
 $$
 \xymatrix{
&K(H)\ar@{-}[dl] \ar@{-}[dr]&\\
K(H_1,\mu_{\ell^{m_2}})\ar@{-}[dr]& & K(H_2)=K(H_2,\mu_{\ell^{m_2}})\ar@{-}[dl] \\
 &K(\mu_{\ell^{m_2}})&}
 $$
Comme l'extension $K(\mu_{\ell^{m_2}})/K(\mu_{\ell^{m_1}})$ est galoisienne et comme $K(H_1)\cap K(\mu_{\ell^{m_2}})=K(\mu_{\ell^{m_1}})$, on en tire que  le diagramme ci-dessous est un parall\'elogramme.
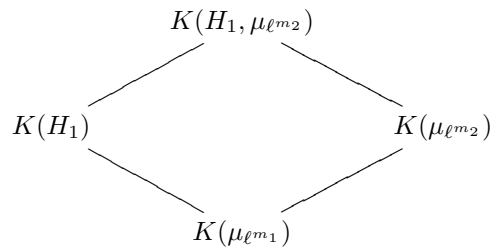
 $$
 \xymatrix{
&K(H_1,\mu_{\ell^{m_2}})\ar@{-}[dl] \ar@{-}[dr]&\\
K(H_1)\ar@{-}[dr]& & K(\mu_{\ell^{m_2}})\ar@{-}[dl] \\
 &K(\mu_{\ell^{m_1}})&}
 $$
D'o\`u l'on tire le parall\'elogramme recherch\'e~:
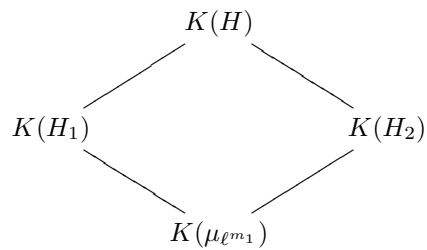
 $$
 \xymatrix{
&K(H)\ar@{-}[dl] \ar@{-}[dr]&\\
K(H_1)\ar@{-}[dr]& & K(H_2)\ar@{-}[dl] \\
 &K(\mu_{\ell^{m_1}})&}
 $$
 Consid\'erons maintenant le diagramme~:
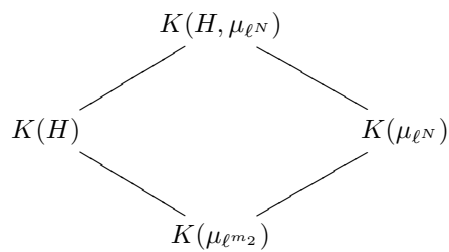
 $$
 \xymatrix{
&K(H,\mu_{\ell^{N}})\ar@{-}[dl] \ar@{-}[dr]&\\
K(H)\ar@{-}[dr]& & K(\mu_{\ell^{N}})\ar@{-}[dl] \\
 &K(\mu_{\ell^{m_2}})&}
 $$
\noindent On peut \'ecrire
\begin{align*}
 [K(H,\mu_{\ell^{N}}):K(\mu_{\ell^{N}})]	& =[K(H_1,\mu_{\ell^{N}}):K(\mu_{\ell^{N}})][K(H_2,\mu_{\ell^{N}}):K(\mu_{\ell^{N}})]\\
											& =[K(H_1):K_1][K(H_2):K_2]=[K(H):K_2]
\end{align*}
\noindent Ainsi le diagramme est un parall\'elogramme et
la deuxi\`eme affirmation de la proposition, \`a savoir, $K(H)\cap K(\mu_{\ell^{\infty}})=K(\mu_{\ell^{m_2}})$ en d\'ecoule donc.
$\findemo$

\medskip

\noindent Une induction ais\'ee permet alors de montrer le th\'eor\`eme suivant~:

\medskip

\begin{theo} Soient $A_1,\dots,A_r$ des vari\'et\'es ab\'eliennes d\'efinies sur $K$ et v\'erifiant la propri\'et\'e $\mu$~:  pour tout $i$ et tout groupe fini $H_i\subset A_i[\ell^{\infty}]$, on a
\[K(H_i)\cap K(\mu_{\ell^{\infty}})\cong K(\mu_{\ell^{m_i}})\ ; \]
\noindent ainsi que l'identit\'e (\ref{hodgehyp}) o\`u $A:=A_1\times\cdots\times A_r$~:
\[\Gal(K(A[\ell^{\infty}])/K(\mu_{\ell^{\infty}}]))\cong \prod_{i=1}^r\Gal(K(A_i[\ell^{\infty}])/K(\mu_{\ell^{\infty}}]))\]
\noindent Alors si $m:=\max m_i $,   tout groupe fini $H=H_1\times\cdots\times H_r\subset A[\ell^{\infty}]$   on a $K(H)\cap K(\mu_{\ell^{\infty}})\cong K(\mu_{\ell^m})$ et
\[[K(H):K(\mu_{\ell^{m}})]\gg\ll\prod_{i=1}^{r}[K(H_i):K(\mu_{\ell^{m_i}})].\]
\end{theo}
\demo Posons $m':=\max\{m_1,\dots,m_{r-1}\}$, de sorte que $m=\max(m',m_r)$. Notons
$A'_1=A_1\times\cdots\times  A_{r-1}$ et $A'_2=A_r$ et enfin $H'_1=H_1\times H_{r_1}$ et $H'_2=H_r$.  Par r\'ecurrence, on sait que $K(H'_1)\cap K(\mu_{\ell^{\infty}})\cong K(\mu_{\ell^{m'}})$ et que
$$[K(H'_1):K(\mu_{\ell^{m'}})]\gg\ll\prod_{i=1}^{r-1}[K(H_i):K(\mu_{\ell^{m_i}})].$$
En appliquant la proposition~\ref{casr=2} \`a $A=A'_1\times A'_2$ et aux sous-groupes $H'_1,H'_2$, on obtient bien que
$[K(H'_1\times H'_2):K(\mu_{\ell^{\infty}})]\cong K(\mu_{\ell^{m}})$ et que
$$[K(H'_1\times H'_2):K(\mu_{\ell^{m}})]\gg\ll [K(H'_1):K(\mu_{\ell^{m'}})][K(H'_2):K(\mu_{\ell^{m_r}})]
\gg\ll \prod_{i=1}^{r}[K(H_i):K(\mu_{\ell^{m_i}})].$$
$\findemo$

\medskip

\noindent R\'e\'enon\c{c}ons la conclusion du th\'eor\`eme sous une forme moins sym\'etrique mais qui est celle que nous utiliserons par le suite.

\begin{cor}\label{c1} Avec les notations et hypoth\`eses du th\'eor\`eme pr\'ec\'edent, et en notant quitte \`a r\'eindexer, $m_r:=\max_{1\leq i\leq r}m_i$, et on a
\[[K(H):K]\gg\ll\prod_{i=1}^r[K(H_i):K]\prod_{i=1}^{r-1}[K(\mu_{\ell^{m_i}}):K]^{-1}.\]
\end{cor}

\subsection{La propri\'et\'e $\mu$ pour les vari\'et\'es ab\'eliennes et $H=A[\ell^{n}]$\label{SLG}}

\noindent L'existence de l'accouplement de Weil implique que,  pour toute vari\'et\'e ab\'elienne $A$ d\'efinie sur $K$, on  $\mu_N\subset K(A[N])$ et en particulier $K(\mu_{\ell^{m}}) \subset K(A[\ell^m])\cap K(\mu_{\ell^{\infty}})$; le but de ce paragraphe est de montrer que l'on a en fait \'egalit\'e, \`a indice fini pr\`es.

\begin{prop} Soient $A$ une vari\'et\'e ab\'elienne d\'efinie sur $K$ , il existe $c_A$ telle que~:
\[\forall \ell \text{ premier },\ \forall m\in\N,\ \  \left[K(A[\ell^m])\cap K(\mu_{\ell^{\infty}}): K(\mu_{\ell^{m}})\right]\leq c_A \]
\end{prop}
\demo
Reprenons la repr\'esentation galoisienne $\rho_{\ell}:\Gal(\overline{K}/K)\rightarrow\GL_{2g}(\Z_{\ell})$ et notons pour abr\'eger $G=G_{\ell}$ l'image, qui est un sous-groupe ferm\'e de $\GL_{2g}(\Z_{\ell})$ et m\^eme, en tenant compte d'une polarisation du groupe des similitudes symplectiques $\GSp_{2g}(\Z_{\ell})$. Ce dernier groupe est muni d'un homomorphisme canonique $\mu:\GSp_{2g}\rightarrow\G_m$ (d\'efini par $\mu(\sigma)e\cdot e'=\sigma(e)\cdot \sigma(e')$). On sait que $\mu\circ\rho_{\ell}$ est le caract\`ere cyclotomique $\chi_{\ell}:\Gal(\overline{K}/K)\rightarrow\Z_{\ell}^{\times}$ dont le noyau correspond \`a l'extension $K(\mu_{\ell^{\infty}})$. Notons $C$ le centre de $\GL_{2g}$, \textit{i.e.} le sous-groupe des homoth\'eties.  D'apr\`es \cite{bogo}, le groupe $G$ contient un sous-groupe d'indice fini dans $C(\Z_{\ell})$ et, d'apr\`es \cite{college8586} cet indice est m\^eme born\'e ind\'ependamment de $\ell$. On a donc
\[ (C(\Z_{\ell}):G\cap C(\Z_{\ell}))\leq c_A.\]
Introduisons les trois sous-groupes $\SLG:=G\cap\SL_{2g}(\Z_{\ell})$, et, pour tout entier $m\geq 1$,
\[H_{\ell^m}:=\Ker\left\{Ê\mu\!\!\!\!\!\mod\ell^m : G\longrightarrow\left(\Z/\ell^m\Z\right)^{\times}ÊÊ\right\}.\]
et 
\[G_{\ell^m}:=\Ker\left\{Ê\red\!\!\!\!\!\mod\ell^m : G\longrightarrow\GL_{2g}(\Z/\ell^m\Z)ÊÊ\right\}.\]
L'extension de $K$ correspondant \`a $H_{\ell^m}$ est $K(\mu_{\ell^m})$ et celle correspondant \`a $G_{\ell^m}$ est $K(A[\ell^m])$.
L'\'enonc\'e de la proposition se traduit alors en l'\'egalit\'e, \`a indice fini born\'e  ind\'ependamment de $\ell$ et $m$:
\[H_{\ell^m}\simeq\langle G_{\ell^m},\SLG\rangle.\]
Consid\'erons l'isog\'enie de groupes alg\'ebriques donn\'ee par le produit,~:
\[\psi: C\times\SL_{2g}\rightarrow\GL_{2g}.\]
D'apr\`es le r\'esultat de Bogomolov compl\'et\'e par Serre rappel\'e ci-dessus, $\psi$ induit un quasi-isomor\-phisme (noyau et conoyau fini born\'e ind\'ependamment de $\ell$) entre $C(\Z_{\ell})\times \SLG$
et $G$. On trouve
\[\psi^{-1}(G_{\ell^m})=\left\{ (\lambda,h)\in \Z_{\ell}^{\times}\times G\;|\; \lambda h\equiv \text{Id}\mod\ell^m  \right\}\]
\[\psi^{-1}(\SLG)=\left\{ (\lambda,h)\in \Z_{\ell}^{\times}\times G\;|\; \lambda^{2g}=1  \right\}\]
\[\psi^{-1}(H_{\ell^m})=\left\{ (\lambda,h)\in \Z_{\ell}^{\times}\times G\;|\; \lambda^{2g}\equiv 1\mod\ell^m  \right\}\]
et il est clair que le dernier groupe est \'egal au produit des deux pr\'ec\'edents, ce qui ach\`eve la d\'emonstration. 
$\findemo$

\section{Preuve du th\'eor\`eme \ref{th1}}

\noindent Soit $A=\prod_{i=1}^nE_i^{n_i}$ un produit de courbes elliptiques sur $\overline{\Q}$ deux \`a deux non-isog\`enes. Nous allons prouver le th\'eor\`eme \ref{th1} annonc\'e dans l'introduction, \`a savoir~:

\[\gamma(A)=\alpha(A).\]

\noindent Rappelons que par ce qui pr\'ec\`ede, il suffit de montrer que 
\[|H|\ll[K(H):K]^{\alpha(A)}\]
\noindent pour $H\subset A[\ell^{\infty}]$ de la forme $\prod_{i=1}^n H_i^{n_i}$ avec $H_i\subset E_i[\ell^{\infty}]$ ; et ce pour tout premier $\ell$.

\medskip

\noindent Nous fixons dans la suite un nombre premier $\ell$.

\subsection{La propri\'et\'e $\mu$ pour les courbes elliptiques}

\noindent Soit $E/K$ une courbe elliptique d\'efinie sur un corps de nombres $K$. Soit $H\subset E[\ell^{\infty}]$ un sous-groupe fini du groupe des points de $\ell^{\infty}$-torsion de $E$ sur $\overline{K}$. Le groupe $H$ est isomorphe au produit $\Z/\ell^{m}\Z\times\Z/\ell^{n}\Z$ avec $0\leq m\leq n$. 

\begin{prop}\label{petitetorsion} La courbe elliptique $E$ v\'erifie la propri\'et\'e $\mu$. Pr\'ecis\'ement, avec les notations pr\'e\-c\'e\-den\-tes, on a
\[\left[K\left(\mu_{\ell^{\infty}}\right)\cap K(H):K\left(\mu_{\ell^{m}}\right)\right]=O(1).\]
et en particulier~:
\[\left[K\left(\mu_{\ell^{\infty}}\right)\cap K(H):K\right]\gg\ll \ell^{m}.\]
\end{prop}


\noindent Notons d\'ej\`a que la deuxi\`eme formule d\'ecoule de la premi\`ere car pour tout entier $r\geq 0$, l'extension $K(\mu_{\ell^r})/K$ est ab\'elienne de degr\'e $\phi(\ell^r)\gg\ll\ell^r$.
Commen\c{c}ons par introduire une d\'efinition dans le cas CM.

\medskip

\defi \label{dep}Soit $E/K$ une courbe elliptique \`a multiplication complexe, et soit $\ell$ un nombre premier. Notons $K_0$ le corps quadratique imaginaire de multiplication complexe. Nous dirons que \textit{$E/K$ est d\'eploy\'ee en $\ell$} si le tore de Mumford-Tate $\MT(E)$ (qui n'est, dans ce cas, autre que $\Res_{K_0/\Q}(\G_{m,K_0})$ restriction des scalaires de $K_0$ \`a $\Q$) est tel que $\MT(E)\times_{\Q} \Q_{\ell}$ est isomorphe au produit $\G_{m,\Q_{\ell}}^2$. Dans le cas contraire nous dirons que \textit{$E/K$ n'est pas d\'eploy\'ee en $\ell$}.

\medskip

\noindent Dans le cas sans multiplication complexe (et dans le cas CM d\'eploy\'e) la preuve de la proposition \ref{petitetorsion} repose sur un petit lemme de groupes. Le cas CM non-d\'eploy\'e se traite diff\'eremment. Commen\c{c}ons donc par ce dernier cas : 

\medskip

\noindent \textbf{Preuve de la proposition \ref{petitetorsion} dans le cas CM non-d\'eploy\'e :} On suppose ici que la courbe elliptique $E$ est de type CM et non-d\'eploy\'ee en $\ell$, au sens de la d\'efinition \ref{dep}. Dans ce cas la repr\'esentation naturelle $\rho_{\Q_{\ell}} : \MT(E)(\Q_{\ell})\rightarrow \Aut(V_{\ell})$ est irr\'eductible. Le paragraphe 7.1 (proposition 7.4) de \cite{ratazziens} nous indique alors que $E[\ell^m]$ est inclus dans $H$ et d'indice fini $c_0$ (ne d\'ependant que de $E/K$ et pas de $H$,$\ell$, $m$) dans $H$. On a donc la s\'erie d'inclusions
\[
K \subset  K\left(\mu_{\ell^m}\right)=K\left(E[\ell^m]\right)\cap K\left(\mu_{\ell^{\infty}}\right)	\subset	 K(H)\cap K\left(\mu_{\ell^{\infty}}\right) \subset K(E[c_0\ell^{m}])\cap K\left(\mu_{\ell^{\infty}}\right)= K\left(\mu_{\ell^{m}}\right)\] 
\noindent La premi\`ere \'egalit\'e ainsi que la derni\`ere \'etant des \'egalit\'es \`a indice fini pr\`es, au sens de la convention du paragraphe \ref{conv}. Elles d\'ecoulent du paragraphe \ref{SLG} indiquant pr\'ecis\'ement que la propri\'et\'e $\mu$ est valable pour $H=A[\ell^m]$ pour toute vari\'et\'e ab\'elienne $A/K$ et d\'eterminant l'exposant $m$ correspondant. \hfill$\Box$

\medskip

\noindent Passons au cas des courbes elliptiques sans multiplication complexe.
\noindent D\'efinissons les sous-groupes suivants de $\GL_2(\Z_{\ell})$.

\[ G_{m,n}:=\left\{\begin{pmatrix}
a	&	b\\
c	&	d\\
\end{pmatrix}\;\text{\Large{$|$}}\; a-1\equiv c\equiv 0 \mod \ell^m,\; b\equiv d-1\equiv 0\mod \ell^{n}\right\}\]

\[\Gamma:=\SL_2(\Z_{\ell})\text{ et } \Gamma_m:=\left\{U\in\GL_2(\Z_{\ell})\;|\;\det(U)\equiv 1\mod \ell^m\right\}.\]
\noindent Nous utiliserons le lemme suivant~:

\begin{lemme}\label{gammamn} Si $m\leq n$, on a l'\'egalit\'e $G_{m,n}\cdot\Gamma=\Gamma_m$.
\end{lemme}
\demo Soit $M=
\begin{pmatrix}
a	&	b	\\
c	&	d	\\
\end{pmatrix}\in\GL_2(\Z_{\ell})$ 
telle que $\det(M)=ad-bc=1\mod \ell^{m}$. La matrice $\begin{pmatrix}
\det(M)^{-1}	&	0	\\
0	&	1	\\
\end{pmatrix}M$ est dans $\SL_2(\Z_{\ell})$.\hfill$\Box$

\medskip

\noindent \textbf{Preuve de la proposition \ref{petitetorsion} dans le cas sans multiplication complexe :} C'est une cons\'equence du lemme pr\'ec\'edent et des r\'esultats connus sur les groupes de Galois des points de $\ell^{\infty}$-torsion des courbes elliptiques sans multiplication complexe. Notons
\[G_1=\Gal\left(K\left(E\left[\ell^{\infty}\right]\right)/K\left(\mu_{\ell^{\infty}}\right)\right)\text{ et }G_H=\Gal\left(K\left(E\left[\ell^{\infty}\right]\right)/K(H)\right).\]
\noindent Commen\c{c}ons par rappeler un r\'esultat de Serre sur la conjecture de Mumford-Tate. Notons $G_{\ell}:=\rho_{\ell}(G_K)$ l'image du groupe de Galois $G_K$ par la repr\'esentation $\ell$-adique $\rho_{\ell} : G_K\rightarrow \Aut(T_{\ell}(E))$ donn\'ee par l'action de Galois sur les points de $E[\ell^{\infty}]$. La courbe $E/K$ \'etant sans multiplication complexe, on sait (d'apr\`es le th\'eor\`eme 3 et son corollaire 1 de \cite{serregalois} p. 299--300) que pour presque tout premier $\ell$, on a 
\[G_{\ell}=\GL_2(\Z_{\ell}).\]
\noindent Ainsi au vu de ce que l'on veut montrer, on peut supposer que $G_{\ell}=\GL_2(\Z_{\ell})$. Le groupe $H$ est de la forme
\[H=\langle P_1\rangle\oplus\langle P_2\rangle\simeq \Z/\ell^{m}\Z\times\Z/\ell^{n}\Z.\]
\noindent Par ailleurs, le groupe de Galois associ\'e, $G_H$, tel que $[G_{\ell}:G_H]=[K(H):K]$ est d\'efini par 
\[G_H=\left\{\sigma\in G_{\ell}\ / \ \sigma_H=\text{Id}_H\right\}.\]
\noindent On se donne une base de $T_{\ell}(E)$, $\{\hat{P_1},\hat{P_2}\}$ telle que $P_1=\hat{P_1}\mod \ell^{m}$ et $P_2=\hat{P_2}\mod \ell^{n}$. On a ainsi l'identification
\[G_H=\left\{\sigma\in \GL_2(\Z_{\ell})\, | \, \sigma \hat{P_1}=\hat{P_1}\mod \ell^{m},\ \text{et}\  \sigma \hat{P_2}=\hat{P_2}\mod \ell^{n}\right\}.\]
\noindent On peut encore r\'e\'ecrire ceci sous la forme
\[G_H=\left\{ \begin{pmatrix}
a	& b	\\
c	& d	\\
\end{pmatrix}\in \GL_2(\Z_{\ell})\ / \ a-1=c=0\mod \ell^{m}\ \text{ et }\ b=d-1=0\mod \ell^{n}\right\}=G_{m,n}.\]
\noindent Par ailleurs, comme le groupe de Mumford-Tate v\'erifie $\MT(E)=\G_m\cdot\Hdg(E)$, on a : $G_1\subset \SL_2(\Z_{\ell})$, et cette inclusion est de plus de conoyau de cardinal major\'e ind\'e\-pen\-damment de $\ell$ (il s'agit du corollaire 2 p. 300 de \cite{serregalois}). Notons $K_1=K(\mu_{\ell^{\infty}})$, $K_2=K(H)$ et $L=K(E[\ell^{\infty}])$. L'extension $L/K$ est galoisienne, donc les extensions $L/K_i$ aussi, de groupe de Galois $G_i$ pour $i\in\{1,2\}$. Ainsi on a
\[\Gal(L/K_1\cap K_2)=G_1\cdot G_H\subset \Gamma_m:=\left\{M\in \GL_2(\Z_{\ell})\ |\ \det M=1\mod\ell^{m}\right\}.\]
\noindent Le lemme \ref{gammamn} pr\'ec\'edent nous assure que cette inclusion est de conoyau fini born\'e ind\'ependamment de $\ell$. On en d\'eduit que
$[K_1\cap K_2:K(\mu_{\ell^m})]$ est born\'e ind\'ependamment de $\ell$. Ceci conclut la preuve dans ce cas.\hfill$\Box$

\medskip

\noindent \textbf{Preuve de la proposition \ref{petitetorsion} dans le cas CM d\'eploy\'e : } il nous reste \`a traiter le cas o\`u $E$ a multiplication complexe et o\`u $\ell$ est tel que le groupe de Mumford-Tate est d\'eploy\'e sur $\Q_{\ell}$, donc isomorphe \`a $\G_{m,\Q_{\ell}}^2$. Dans ce cas, le th\'eor\`eme \ref{trr} indique que $G_{\ell}$ est ouvert, d'indice fini dans $\Z_{\ell}^{\times}\times\Z_{\ell}^{\times}$ plong\'e diagonalement dans $\GL_2(\Z_{\ell})$. Le groupe $G_H$ est donc (\`a indice fini pr\`es) le groupe
\[\left\{ \begin{pmatrix}
a	& 0	\\
0	& d	\\
\end{pmatrix}\in \GL_2(\Z_{\ell})\ | \ a=1\mod \ell^{m}\ \text{ et }\ d=1\mod \ell^{n}\right\}.\]
\noindent Avec les m\^emes notations que ci-dessus, on a $G_1\subset \left\{ \begin{pmatrix}
a	& 0	\\
0	& a^{-1}	\\
\end{pmatrix}\in \GL_2(\Z_{\ell})\ | \ a\in\Z_{\ell}^{\times}\right\}$, et cette inclusion est de conoyau de cardinal major\'e ind\'e\-pen\-damment de $\ell$. De plus, 
\[\Gal(L/K_1\cap K_2)=G_1\cdot G_H\subset \Gamma_m:=\left\{ \begin{pmatrix}
a	& 0	\\
0	& d	\\
\end{pmatrix}\in \GL_2(\Z_{\ell})\ |\ ad=1\mod\ell^{m}\right\}.\] 
\noindent L'\'equivalent du lemme de groupes pr\'ec\'edent est ici le suivant : si $M =\begin{pmatrix}
a	& 0	\\
0	& d	\\
\end{pmatrix}$ est telle que $ad=1\mod\ell^m$, alors en la multipliant par la matrice $\begin{pmatrix}
a^{-1}	& 0	\\
0	& a	\\
\end{pmatrix}$ on obtient la matrice $\begin{pmatrix}
1	& 0	\\
0	& ad	\\
\end{pmatrix}$ avec $ad=1\mod\ell^m$. Ceci montre que $[K(\mu_{\ell^{\infty}})\cap K(H):K(\mu_{\ell^m})]$ est born\'e ind\'ependamment de $\ell$.

\noindent Ceci conclut la preuve.\hfill$\Box$

\subsection{Combinatoire}

\noindent Revenons \`a notre situation initiale : soit $A=\prod_{i=1}^nE_i^{n_i}$ produit de courbes elliptiques, les $E_i$ \'etant deux \`a deux non-isog\`enes, et soit pour tout $i$ des sous-groupes $H_i\subset E_i[\ell^{\infty}]$ finis. Rappelons que par le th\'eor\`eme \ref{MTfort} et le lemme \ref{produit}, on sait que $A_{\text{red}}:=\prod_{i=1}^nE_i$ v\'erifie les hypoth\`eses de la proposition \ref{h4}. Donc $A_{\text{red}}$ et $A$ v\'erifient \'egalement la conclusion de la proposition \ref{h4}. De plus par la proposition \ref{petitetorsion}, chacune des $E_i$ v\'erifie la propri\'et\'e $\mu$. Finalement ceci prouve que la vari\'et\'e ab\'elienne $A$ est justiciable du corollaire \ref{c1} que nous utiliserons librement dans la suite.

\medskip

\noindent Il conviendra de distinguer selon que la courbe elliptique $E_i$ est ou non de type CM. Quitte \`a renum\'eroter les $E_i$, on peut \'ecrire  $A$ sous la forme 
\[A=B\times C\ \text{ avec }\ B=\prod_{i=1}^m B_i^{v_i} \text{ et } C=\prod_{i=1}^nC_i^{u_i}.\]
\noindent Les $B_i$ sont des courbes elliptiques avec multiplication complexe et les $C_i$ sont des courbes elliptiques sans multiplication complexe. Ces courbes sont deux \`a deux non-isog\`enes. Les $u_i$ et $v_i$ sont dans $\N-\{0\}$ et $n$ et $m$ sont des entiers positifs ou nuls. On se donne un groupe $H\subset A[\ell^{\infty}]$ de la forme
\[H=H_B\times H_C, \text{ o\`u } H_B=\prod_{i=1}^m H_{B,i}^{v_i}, \text{ et o\`u }H_C=\prod_{i=1}^nH_{C,i}^{u_i}\]
\noindent avec, quitte \`a renum\'eroter,
\[\forall i\in\{1,\ldots,n\}\ \ H_{C,i}\simeq\Z/\ell^{c_i}\Z\times \Z/\ell^{c_{i+n}}\Z\ \text{ et }\ 0\leq c_i\leq c_{i+n}\ \text{ et }\ c_1\leq \ldots\leq c_n,\]
\noindent et 
\[\forall i\in\{1,\ldots,m\}\ \ H_{B,i}\simeq\Z/\ell^{b_i}\Z\times \Z/\ell^{b_{i+m}}\Z\ \text{ et }\ 0\leq b_i\leq b_{i+n}\ \text{ et }\ b_1\leq \ldots\leq b_m.\]

\noindent En notant $\log_{\ell}$ le logarithme en base $\ell$, on a donc
\begin{equation}\label{H}
\log_{\ell}|H|=\sum_{i=1}^n(c_i+c_{i+n})u_i+\sum_{i=1}^m(b_i+b_{i+m})v_i.
\end{equation}

\noindent Par ailleurs la proposition \ref{petitetorsion} dit alors que pour tout $i\in\{1,\dots,m\}$ (resp. tout $j\in\{1,\dots,n\}$), on a $K(H_{B,i})\cap K\left(\mu_{\ell^{\infty}}\right)=K\left(\mu_{\ell^{b_i}}\right)$ \`a indice fini pr\`es (resp. $K(H_{C,j})\cap K\left(\mu_{\ell^{\infty}}\right)=K\left(\mu_{\ell^{c_j}}\right)$) et on a donc des quasi-\'egalit\'es du type~:

\[\left[K(H_{B,i})\cap K\left(\mu_{\ell^{\infty}}\right):K\right]\gg\ll\left[K(\mu_{\ell^{b_i}}):K\right]\gg\ll\ell^{b_i}.\]

\noindent Introduisons de plus la notation 
\[\beta=\min\{c_n,b_m\}\ \text{ si $m$ et $n$ sont non nuls, et }\ \ \ \beta=0 \text{ si $m$ ou $n$ est nul.}\]

\noindent En appliquant le corollaire \ref{c1} nous obtenons ainsi
\begin{equation}\label{fait}
[K(H):K]\gg\ll\left(\prod_{i=1}^m[K(H_{B,i}):K]\prod_{j=1}^n[K(H_{C,j}):K]\right)\ell^{-\sum_{i=1}^{m-1}b_i-\sum_{j=1}^{n-1}c_j-\beta}.
\end{equation}

\noindent Par ailleurs, en utilisant les r\'esultats de \cite{ratazziens} pour les courbes elliptiques sans CM et pour les courbes elliptiques avec CM, on a :
\begin{enumerate}
\item pour tout $i$, $\log_{\ell}[K(H_{C,i}):K]\geq 2(c_i+c_{i+n})+ O(1)$.
\item pour tout $i$, $\log_{\ell}[K(H_{B,i}):K]\geq b_i+b_{i+n}+ O(1)$ (dans \cite{ratazziens} ceci est prouv\'e si $H_{B,i}$ est galoisien mais on a vu dans la r\'eduction du paragraphe \ref{redgalois} que ceci entra\^ine l'in\'egalit\'e dans le cas g\'en\'eral).
\end{enumerate}
\noindent En utilisant ceci, ainsi que l'encadrement (\ref{fait}) on obtient
\begin{align*}
\log_{\ell}\left([K(H):K]\right)	& =\sum_{i=1}^n\log_{\ell}[K(H_{C,i}):K]+\sum_{i=1}^m\log_{\ell}[K(H_{B,i}):K]-\sum_{i=1}^{n-1}c_i-\sum_{i=1}^{m-1}b_i-\beta+O(1)\\
									& \geq 2\sum_{i=1}^{n}(c_i+c_{i+n})-\sum_{i=1}^{n-1}c_i+\sum_{i=1}^{m}(b_i+b_{i+m})-\sum_{i=1}^{m-1}b_i-\beta+O(1)\\
									& \geq c_n+\sum_{i=1}^nc_i+2\sum_{i=1}^nc_{i+n}+\sum_{i=1}^mb_{i+m}+b_m-\beta+O(1).
\end{align*}
\noindent Introduisons la notation
\[m(A):=\max\left\{\frac{\sum_{i=1}^n(c_i+c_{i+n})u_i+\sum_{i=1}^m(b_i+b_{i+m})v_i}{c_n+\sum_{i=1}^nc_i+2\sum_{i=1}^nc_{i+n}+\sum_{i=1}^mb_{i+m}+b_m-\beta}\right\},\]
\noindent le max portant sur les conditions $\ 0\leq b_i\leq b_{i+m},\ \ 0\leq c_i\leq c_{i+n},\ \ b_1\leq\ldots\leq b_m,\ \ c_1\leq\ldots\leq c_n$.

\medskip

\noindent Finalement la preuve du th\'eor\`eme \ref{th1} se ram\`ene \`a la proposition combinatoire suivante~:

\begin{prop}\label{combin} $\alpha(A)\geq m(A)$.
\end{prop}
\demo Nous allons distinguer 3 cas, selon que $n=0$ ou $m=0$ ou $n,m\geq 1$.

\medskip

\noindent 1. Commen\c{c}ons pas le cas o\`u l'on suppose $m=0$, \textit{i.e.} $A$ est un produit de courbes elliptiques sans CM (et donc $\beta=0$). On a
\[\alpha(A)=\max_{\emptyset\not=I\subset\{1,\ldots, n\}}\frac{2\sum_{i\in I}c_i}{\dim \MT\left(\prod_{i\in I}E_i^{u_i}\right)}.\]
\noindent Or par les rappels sur le groupe de Mumford-Tate (lemmes \ref{puissance0} et \ref{produit}), on a
\[\dim \MT(\prod_{i\in I}E_i^{u_i})=\dim\MT(\prod_{i\in I}E_i)=1+\dim\Hdg(\prod_{i\in I}E_i)=1+\sum_{i\in I}\dim \Hdg(E_i)=1+3|I|.\]
\noindent (La derni\`ere \'egalit\'e vient de ce que $\MT(A)=\GL_{2}$ dans le cas sans multiplication complexe, donc $\dim\Hdg(A)=3$). Finalement on a 
\begin{equation}\label{e0}
\alpha(A)=\max_{\emptyset\not=I\subset\{1,\ldots, n\}}\frac{2\sum_{i\in I}u_i}{1+3|I|}\geq \max_{1\leq i\leq n}\frac{u_i}{2}.
\end{equation}
\noindent Nous pouvons maintenant passer \`a la preuve de l'in\'egalit\'e $\alpha(A)\geq m(A)$ proprement dite. Celle-ci est \'equivalente \`a montrer que
\begin{equation}\label{e1}
\sum_{i=1}^{n-1}c_i(u_i-\alpha(A))\leq c_n(2\alpha(A)-u_n)+\sum_{i=1}^nc_{i+n}(2\alpha(A)-u_i).
\end{equation}
\noindent Or on a 
\[\sum_{i=1}^{n-1}c_i(u_i-\alpha(A))\leq\sum_{i\in I}c_i(u_i-\alpha(A))\  \text{ avec }\ I=\{i\leq n-1\ |\ u_i-\alpha(A)\geq 0\}.\]
\noindent De plus, comme $2\alpha(A)-u_i\geq 0$ et $c_{i+n}\geq c_i$, on a
\begin{align*} 
c_n(2\alpha(A)-u_n)+\sum_{i=1}^nc_{i+n}(2\alpha(A)-u_i)	& \geq c_n(4\alpha(A)-2u_n)+\sum_{i=1}^{n-1}c_{i}(2\alpha(A)-u_i)\\
														& \geq c_n(4\alpha(A)-2u_n)+\sum_{i\in I}c_{i}(2\alpha(A)-u_i).
\end{align*}
\noindent Donc pour prouver (\ref{e1}), il suffit de montrer que 
\[\sum_{i\in I}c_i(u_i-\alpha(A))\leq c_n(4\alpha(A)-2u_n)+\sum_{i\in I}c_{i}(2\alpha(A)-u_i),\]
\noindent autrement dit que
\[\sum_{i\in I}c_i(2u_i-3\alpha(A))\leq c_n(4\alpha(A)-2u_n).\]
\noindent Mais on a \'egalement
\[\sum_{i\in I}c_i(2u_i-3\alpha(A))\leq \sum_{i\in I'}c_i(2u_i-3\alpha(A))\ \text{ o\`u }\ I'=\{i\in I\ | 2u_i-3\alpha(A)>0\}.\]
\noindent Ainsi il suffit, pour prouver (\ref{e1}), de montrer que
\[\sum_{i\in I'}c_i(2u_i-3\alpha(A))\leq c_n(4\alpha(A)-2u_n).\]
\noindent Mais $\sum_{i\in I'}c_i(2u_i-3\alpha(A))\leq \sum_{i\in I'}c_n(2u_i-3\alpha(A))$, donc il nous reste \`a montrer que 
\[2\sum_{i\in I'\cup\{n\}}u_i\leq 4\alpha(A)+3|I'|\alpha(A)=\left(3|I'\cup\{n\}|+1\right)\alpha(A).\]
\noindent Ceci est vrai comme le montre (\ref{e0}) avec l'ensemble $I'\cup\{n\}$.

\medskip

\noindent 2.  Supposons maintenant que $n=0$, \textit{i.e.} que $A$ est un produit de courbes elliptiques ayant multiplication complexe (et donc avec $\beta=0$). On a
\[\alpha(A)=\max_{\emptyset\not=I\subset\{1,\ldots, m\}}\frac{2\sum_{i\in I}b_i}{\dim \MT\left(\prod_{i\in I}E_i^{v_i}\right)}.\]
\noindent Le m\^eme calcul qu'\`a l'\'etape 1. (en utilisant que pour une courbe elliptique CM, le groupe $\MT(E)$ est de dimension $2$) donne cette fois
\begin{equation}\label{e2}
\alpha(A)=\max_{\emptyset\not=I\subset\{1,\ldots, m\}}\frac{2\sum_{i\in I}v_i}{1+|I|}\geq \max_{1\leq i\leq m}v_i.
\end{equation}
\noindent Nous pouvons maintenant passer \`a la preuve de l'in\'egalit\'e $\alpha(A)\geq m(A)$ proprement dite. Celle-ci est \'equivalente \`a montrer que
\[\sum_{i=1}^{m}(b_i+b_{i+m})v_i\leq \alpha(A)b_m+\alpha(A)\sum_{i=1}^mb_{i+m}\]
\noindent autrement dit \`a 
\begin{equation}\label{e3}
\sum_{i=1}^{m}b_iv_i\leq \alpha(A)b_m+\sum_{i=1}^mb_{i+m}(\alpha(A)-v_i).
\end{equation}
\noindent Comme $b_{i+m}\geq b_i$ et comme $\alpha(A)-v_i\geq 0$ on a $\sum_{i=1}^mb_{i+m}(\alpha(A)-v_i)\geq \sum_{i=1}^mb_{i}(\alpha(A)-v_i)$, donc pour prouver (\ref{e3}), il suffit de montrer que 
\[\sum_{i=1}^{m-1}b_i(2v_i-\alpha(A))\leq 2(\alpha(A)-v_m)b_m.\]
\noindent En introduisant comme auparavant l'ensemble $I=\{i\leq m-1\ |\ 2v_i-\alpha(A)\geq 0\}$, et en utilisant que $b_i\leq b_m$, on voit qu'il suffit de montrer que 
\[2v_m+2\sum_{\in I}v_i\leq\alpha(A)\left(|I|+2\right).\]
\noindent En consid\'erant l'ensemble $I\cup\{m\}$ on constate avec (\ref{e2}) que cette derni\`ere relation est vraie.

\medskip

\noindent 3. Il nous reste maintenant le cas o\`u $m$ et $n$ sont non-nuls. L\`a encore c'est le m\^eme type de calcul qui va permettre de conclure. Cette fois on a la formule 
\begin{equation}\label{e4}
\alpha(A)=\max_{I,J,\ \emptyset\not=I\cup J}\frac{2\sum_{i\in I}u_i+2\sum_{j\in J}v_j}{3|I|+|J|+1}\geq \max\{\max_{1\leq i\leq n}\frac{u_i}{2} ; \max_{1\leq i\leq m}{v_i}\}.
\end{equation}
\noindent L'in\'egalit\'e \`a d\'emontrer est la suivante~:
\begin{equation}\label{e5}
\sum_{i=1}^n(c_i+c_{i+n})u_i+\sum_{i=1}^m(b_i+b_{i+m})v_i\leq \alpha(A)\left[\sum_{i=1}^nc_i+2\sum_{i=1}^nc_{i+n}+\sum_{i=1}^mb_{i+m}+(c_n+b_m-\beta )\right].
\end{equation}
\noindent Comme pr\'ec\'edemment, en utilisant que $\alpha(A)-v_i\geq 0$ et que $b_{i+m}\geq b_i$, on voit qu'il suffit de montrer que
\[\sum_{i=1}^nc_iu_i+\sum_{i=1}^mb_iv_i\leq \alpha(A)\sum_{i=1}^nc_i+(2\alpha(A)-u_i)\sum_{i=1}^nc_{i+n}+(\alpha(A)-v_i)\sum_{i=1}^mb_i+\alpha(A)(c_n+b_m-\beta).\]
\noindent De m\^eme, en utilisant que $2\alpha(A)-u_i\geq 0$ et que $c_{i+n}\geq c_i$, il suffit de montrer que 
\[\sum_{i=1}^{n-1}(2u_i-3\alpha(A))c_i+\sum_{j=1}^{m-1}(2v_j-\alpha(A))b_j\leq c_n(4\alpha(A)-2u_n)+b_m(2\alpha(A)-2v_m)-\beta\alpha(A).\]
\noindent Introduisons les ensembles $I'=\{i\leq n-1\ |\ 2u_i-3\alpha(A)\geq 0\}$ et $J'=\{j\leq m-1\ |\ 2v_j-\alpha(A)\geq 0\}$. Il suffit, pour prouver (\ref{e5}), de montrer que 
\[\sum_{i\in I'}(2u_i-3\alpha(A))c_n+\sum_{j\in J'}(2(v_j-\alpha(A))b_m\leq c_n(4\alpha(A)-2u_n)+b_m(2\alpha(A)-2v_m)-\beta\alpha(A).\]
\noindent En posant $I=I'\cup\{n\}$ et $J=J'\cup\{m\}$, ceci revient \`a prouver que 
\begin{align*}
\left(\sum_{i\in I}2u_i\right)c_n+\left(\sum_{j\in J}2v_j\right)b_m	&	\leq (4\alpha(A)+3\alpha(A)|I'|)c_n+(2\alpha(A)+\alpha(A)|J'|)b_m-\alpha(A)\beta\\
																	&	\leq \alpha(A)\left((3|I|+1)c_n+(|J|+1)b_m-\beta\right).
\end{align*}
\noindent Si $\beta:=\min\{b_m,c_n\}=b_m$, \textit{i.e.} si $b_m\leq c_n$, alors il nous reste \`a montrer que
\[\left(\sum_{i\in I}2u_i\right)c_n+\left(\sum_{j\in J}2v_j\right)b_m\leq \alpha(A)\left[(3|I|+1)c_n+|J|)b_m\right].\]
\noindent Ceci peut se r\'e\'ecrire sous la forme
\[\left(\sum_{j\in J}2v_j-|J|\alpha(A)\right)b_m\leq \left((3|I|+1)\alpha(A)-2\sum_{i\in I}u_i\right)c_n.\]
\noindent Comme $(3|I|+1)\alpha(A)-\sum_{i\in I}2u_i\geq 0$ et comme $b_m\leq c_n$, on voit qu'il suffit de montrer que 
\[\left(\sum_{j\in J}2v_j-|J|\alpha(A)\right)b_m\leq \left((3|I|+1)\alpha(A)-2\sum_{i\in I}u_i\right)b_m.\]
\noindent ceci est entrain\'e par
\[2\sum_{j\in J}v_j+2\sum_{i\in I}u_i\leq (|J|+3|I|+1)\alpha(A).\]
\noindent Or l'in\'egalit\'e (\ref{e4}) nous dit pr\'ecis\'ement que ceci est vrai.

\medskip 

\noindent Le cas o\`u $\beta=c_n$ se traite exactement de la m\^eme fa\c{c}on. Ceci termine la preuve de la proposition \ref{combin} et donc du th\'eor\`eme \ref{th1}.\hfill$\Box$

\vspace{1cm}

\noindent \textbf{Adresse :} Hindry Marc,\\ Institut de Math\'ematiques de Jussieu\\Universit\'e Paris 7 Denis Diderot\\ Case Postale 7012\\ 2, place Jussieu\\ F-75251 PARIS CEDEX 05

\medskip

\noindent Ratazzi Nicolas, \\ Universit\'e Paris-Sud\\ D\'epartement de math\'ematiques, B\^atiment 425\\ 91405 Orsay Cedex \\ FRANCE

\end{document}